\documentclass{amsart}
\usepackage[american]{babel}
\usepackage{amsopn}
\usepackage{amsmath,amssymb}
\usepackage{graphicx}
\usepackage{float}
\usepackage{fancyhdr}
\numberwithin{equation}{section}
\def\R{{\mathbb R}}
\def\E{{\mathbb E}}
\def\N{{\mathbb N}}

\def\spn{{\rm span}}
\def\F{{\mathcal F}}

\def\p{{\partial}}
\def\xL{{\rm L}}
\def\xW{{\rm W}}
\def\xCn{{\rm C}}
\def\s{{\sigma}}
\def\xdif{{\rm d}}
\def\xLtwo{{{\rm L}^2}}
\def\xHone{{{\rm H}^1}}

\def\L2{{\rm L}^2(\R^d)}
\def\xLn{{\rm L}}

\def\xHn{{\rm H}}
\def\xtr{{\rm tr}}
\def\H1{{\xHone(\R^d)}}
\def\spn{{\rm span}}
\def\indic{{\rm {\large 1}\hspace{-2.3pt}{\large
l}}}
\def\xLinfty{{\rm L}^{\infty}}
\def\xim{{\rm im}}
\def\xker{{\rm ker}}
\theoremstyle{plain}
\newtheorem{thrm}{Theorem}[section]
\newtheorem{lmm}[thrm]{Lemma}

\newtheorem{prpstn}[thrm]{Proposition}
\theoremstyle{definition}

\newtheorem{rmrk}[thrm]{Remark}

\title{Uniform large deviations for the nonlinear
Schr\"odinger equation with multiplicative noise}
\author{Eric GAUTIER$^{1,2}$}

\begin{document}

\maketitle  \pagestyle{myheadings} \markboth{{E. Gautier}}{{Uniform
large deviations for the nonlinear Schr\"odinger equation with
multiplicative noise}}
\begin{abstract} Uniform large deviations
for the laws of the paths of the solutions of the stochastic
nonlinear Schr\"odinger equation when the noise converges to zero
are presented. The noise is a real multiplicative Gaussian noise. It
is white in time and colored in space. The path space considered
allows blow-up and is endowed with a topology analogue to a
projective limit topology. Thus a large variety of large deviation
principle may be deduced by contraction. As a consequence,
asymptotics of the tails of the law of the blow-up time when the
noise converges to zero are obtained.\vspace{0.3cm}

\end{abstract}

\footnotetext[1]{CREST-INSEE, URA D2200, 3 avenue Pierre Larousse,
92240 Malakoff, France} \footnotetext[2]{IRMAR, UMR 6625,
Universit\'e de Rennes 1, Campus de Beaulieu, 35042 Rennes cedex,
France; \email{eric.gautier@bretagne.ens-cachan.fr}\vspace{0.2cm}

}
\thispagestyle{empty}

\section{Introduction}\label{s1}
In the present article, the stochastic nonlinear Schr\"odinger (NLS)
equation with a power law nonlinearity and a multiplicative noise is
studied. The deterministic equation occurs as a generic model in
many areas of physics: solid state physics, optics, hydrodynamics,
plasma physics, molecular biology, field theory. It describes the
propagation of slowly varying envelopes of a wave packet in media
with both nonlinear and dispersive responses. In some physical
applications, see \cite{BCIR1,BCIR2,BCGR,DC}, spatial and temporal
fluctuations of the parameters of the medium have to be considered
and sometimes the only source of phase fluctuation that has
significant effect on the dynamics enters linearly as a random
potential. This interaction term may for example account for thermal
fluctuations or inhomogeneities in the medium. The evolution
equation may be written
\begin{equation}
\label{e1} i \frac{\partial u}{\partial t} - (\Delta u
+\lambda|u|^{2\s}u)=u\xi, \quad \lambda=\pm 1,
\end{equation}
where $\xi$ is a real valued Gaussian noise, it is ideally a
space-time white noise with correlation
$$\E\left[\xi(t_1,x_1)\xi(t_2,x_2)\right]=D
\delta_{t_1-t_2}\otimes\delta_{x_1-x_2},$$ $D$ is the noise
intensity, $\delta$ denotes the Dirac mass, $u$ is a complex-valued
function of time and space and the space variables $x$ is in $\R^d$.
When $\lambda=1$ the nonlinearity is called focusing or attractive,
otherwise it is defocusing or repulsive.\vspace{0.3cm}

Since the unbounded operator $-i\Delta$ on $\xHone(\R^d)$ with
domain $\xHn^{3}(\R^d)$ is skew-adjoint, it generates a unitary
group $\left(U(t)=e^{-it\Delta}\right)_{t\in\R}$. Thus
$\left(U(t)\right)_{t\in\R}$ is an isometry on $\xHone(\R^d)$ and
there is no smoothing effect in the Sobolev spaces. We are thus
unable to treat the space-time white noise and will consider a real
valued centered Gaussian noise, white in time and colored in space.
Also when $d>1$, which is for example the case in \cite{BCIR2}, the
nonlinearity is never Lipschitz on the bounded sets of
$\xHone(\R^d)$. Thankfully, the Strichartz estimates, presented in
the next section, show that some integrability property is obtained
through "convolution" with the group. This allows us to treat the
nonlinearity.\vspace{0.3cm}

The well posedness of the Cauchy problem associated to the
deterministic NLS equation depends on the size of $\s$. If
$\s<\frac{2}{d}$, the nonlinearity is subcritical and the Cauchy
problem is globally well posed in $\L2$ or $\H1$. If
$\s=\frac{2}{d}$, critical nonlinearity, or
$\frac{2}{d}<\s<\frac{2}{d-2}$ when $d\geq 3$ or simply
$\s>\frac{2}{d}$ otherwise, supercritical nonlinearity, the Cauchy
problem is locally well posed in $\H1$, see \cite{kato}. In this
latter case, if the nonlinearity is defocusing, the solution is
global. In the focusing case some initial data yield global
solutions while it is known that other initial data yield solutions
which blow up in finite time, see \cite{CAZ,SS}. Two quantities are
invariant on every time interval included in the existence time
interval of the maximal solution. They are the momentum
$$M(u_d^{u_0}(t))=\|u_d^{u_0}(t)\|_{\L2}$$ and the Hamiltonian
$$H(u_d^{u_0}(t))=\frac{1}{2}\int_{\mathbb{R}^d}|\nabla u_d^{u_0}(t)|^2
dx-\frac{\lambda}{2\sigma+2}\int_{\mathbb{R}^d}|u_d^{u_0}(t)|^{2\sigma+2}dx,$$
where we denote by $u_d^{u_0}$ the solution of the deterministic
equation.\vspace{0.3cm}

In this article, we adopt the formalism of stochastic evolution
equations in M-type 2 Banach spaces as presented in \cite{B1,B2},
see also \cite{DPZ} for the case of an Hilbert space. The real
Gaussian noise that we consider hereafter is defined as the time
derivative in the sense of distributions of a Wiener process
$(W(t))_{t\in[0,+\infty)}$ defined on a real separable Banach space.
It is correlated in the space variables and is such that the product
make sense in the space considered for the fixed point argument. We
write $W=\Phi W_c$, where $W_c$ is a cylindrical Wiener process and
$\Phi$ is a particular operator. With the It\^o notations, the
stochastic evolution equation is written
\begin{equation}
\label{e2} i \xdif u -(\Delta u + \lambda |u|^{2\s}u)\xdif t = u
\circ \xdif W.
\end{equation}
The symbol $\circ$ stands for the Stratanovich product. It follows
that the momentum remains an invariant quantity of the stochastic
functional flow. We will use the equivalent It\^{o} form
\begin{equation}
\label{e2} i \xdif u -\left(\Delta u + \lambda
|u|^{2\s}u-\frac{i}{2}uF_{\Phi}\right)\xdif t = u \xdif W,
\end{equation}
where $$F_{\Phi}(x)=\sum_{j\in\N}\left(\Phi e_j
(x)\right)^2,\hspace{0.3cm} x\in\R^d.$$ The initial datum $u_0$ is a
function in $\H1$. We consider solutions of NLS that are weak
solutions in the sense of the analysis of partial differential
equations or equivalently mild solutions which satisfy
\begin{align*}
u^{u_0}(t)=&U(t)u_0-\int_0^tU(t-s)\left(i\lambda|u^{u_0}(s)|
^{2\sigma}u^{u_0}(s)+\frac{1}{2}u^{u_0}(s)F_{\Phi}\right)ds\\
&-i\int_0^tU(t-s)(u^{u_0}(s)dW(s)).
\end{align*}
In \cite{dBD2,dBD1}, local existence and uniqueness for paths in
$\H1$, respectively in  $\xLtwo(\R^d)$, is proved in the stochastic
case when the noise is multiplicative. Global existence when the
nonlinearity is defocusing and in the subcritical case when the
nonlinearity is focusing are also obtained. Investigations on the
blow-up in the stochastic case when the noise is multiplicative
appeared in a series of numerical, see \cite{BSDDM,DMD}, and
theoretical papers, see \cite{dBD3}. The small noise asymptotics is
studied in \cite{EG} in the case of an additive noise. There results
are applied to the study of the blow-up and to the computation of
error probabilities in soliton transmission in fibers. In the second
case solitons are carriers for bits, $d=\sigma=\lambda=1$ and the
$x$ and $t$ variables stand for respectively time and space. To make
a decision on the received pulse, the momentum is computed over a
window. Noise may cause the loss of the signal by essentially
shifting the arrival time (timing jitter) of the soliton and
degrading its momentum or create from nothing a structure that is
large enough and may be mistaken as a soliton. The sample path large
deviation principle (LDP) allowed to obtain similar results on the
tails of the momentum of the pulse as in \cite{FKLT,FKLMT} where the
heuristic method of collective coordinates and the instanton
formalism have been used.\vspace{0.3cm}

In this article, we are again interested in the tails of the laws of
the paths of the solution of this random perturbation of a
Hamiltonian system when the noise goes to zero. Note that in
fiber-optical communications the multiplicative noise may account
for Raman noise, it particularly affects the propagation of short
pulses, see \cite{DC} for more details. However, contrary to the
case of an additive noise, here the momentum of a pulse is
conserved. Nonetheless we could consider the tails of the law of the
timing jitter. We plan to study this application in future works. We
are also able to deduce from the sample path large deviations
results on the asymptotics of the blow-up time, results are stated
in the last section of this article. There also remains many
interesting problems which may be studied when a uniform sample path
large deviation principle is available. Uniformity with respect to
initial data in compact sets is needed when we study the most
probable exit points from the boundary of a more geometric domain
$D$ than the preceding window. The domain may be a compact
neighborhood $D$ of an asymptotically stable equilibrium point. The
corresponding asymptotics of the exit time from the domain may also
be studied, see \cite{DZ,FW} in the case of diffusions and
\cite{CM,Fr} in the case of particular SPDEs. Some results on
asymptotic stability of solitary waves are available for the
translation invariant NLS equation studied herein or with an
additional potential accounting for a defect in the homogeneous
background, see for example \cite{BS,CU,FGJS,SW}. The solution
decomposes then into a solitary wave with temporal fluctuation of
the parameters called modulations and a radiative decaying part.
Contrary to \cite{CM,Fr}, compactness is a real issue in the NLS
equation on an unbounded domain of $\R^d$ when considering random
perturbations since the group does not have a smoothing effect. We
plan to investigate these types of problems in future works. Uniform
LDPs also yield LDPs for the family of invariant measures of Markov
transition semi-groups defined by SPDEs, see \cite{CR2,Sow2} for the
case of some reaction-diffusion equations, when the noise goes to
zero and when the measures weakly converge to a Dirac mass on 0, the
only stationary solution of the deterministic equation. In the case
of the NLS equation there are several invariant measures and little
is known about invariant measures in the stochastic case without
adding terms accounting for damping or viscosity and considering the
space variables in a bounded domain of $\R^d$.\vspace{0.3cm}

The author is aware of three different types of proof for a LDP when
the noise is multiplicative. One uses an extension of Varadhan's
contraction principle that may be found in \cite{DS1} or in
\cite{DZ}[Proposition 4.2.3], see also \cite{FZ}. It requires a
sequence of approximation of the measurable It\^{o} map by
continuous maps uniformly converging on the level sets of the
initial good rate function and that the resulting sequence of family
of laws are exponentially good approximations of the laws of the
solutions. It is needed that the It\^{o} map takes its values in a
metric space and this is not the case when we consider the spaces of
exploding paths, see \cite{AZ1,EG}. The second type of proof is
strongly based on the estimate given in Proposition \ref{prpstn5}
below, it is the one we adopt in the following. We are indebted to
\cite{AZ1} where the proof is written for diffusions that may blow
up in finite time. Some aspects of the proof has been simplified in
\cite{Pr}. Also, some assumptions on the coefficient in front of the
noise were relaxed. Nonetheless the locally Lipschitz conditions on
the drift has been replaced by a globally Lipschitz assumption and
boundedness of the drift and coefficient in front of the noise and
the framework no longer allows SDEs that blow-up in finite time. The
proofs of LDP for stochastic PDEs that we may found for example in
\cite{CW,CR,CM,P1} follow this second type of proof. Note that in
\cite{CW,CM} the approach to the SPDE is based on martingale
measures and the Brownian sheet field whereas in \cite{CR,P1} the
approach has an infinite dimensional flavor, like ours. A third type
of proof has appeared in \cite{LeQiZh} for diffusions. It requires
$\xCn^3$ with linear growth diffusion coefficients and is based on
the continuity theorem of T. Lyons for rough paths in the
$p-$variation topology. Note that the continuity theorem also yields
proofs of support theorems for diffusions. Reference \cite{LGT}
presents in this setting the existence of mild solutions of SPDEs
when the PDE is linear and the semi-group is analytic, which covers
the case of the linear Heat Equation, and when it is driven by a
multiplicative $\gamma$-h\"older time continuous path with value in
a space of distributions. The author is not aware of a proof of a
LDP in the SPDE setting via the Rough Paths theory. Uniform LDPs are
for example proved in \cite{AZ1,DZ,FW} for diffusions, and in
\cite{CM,Sow} for particular SPDEs in spaces of H\"older continuous
functions in both time and space when the space variables are in a
bounded domain.\\

The paper is organized as follows. Section \ref{s2} is devoted to
notations and preliminaries and states the uniform LDP in a space of
exploding paths where blow-up in finite time is allowed. Since the
stronger the topology, the sharper are the estimates, we take
advantage of the variety of spaces where the fixed point can be
conducted, as it has been done in \cite{EG}, due to the
integrability property and endow the space with a topology analogous
to a projective limit topology. The result can be transferred to
weaker topologies or more generally by any family of equicontinuous
mappings using \cite{Sow}[Proposition 5] which is an extension of
Varadhan's contraction principle. The rate function is the infimum
of the $\xLtwo$-norm of the controls, of a deterministic control
problem, producing the prescribed path. The mild solution of the
control problem is called the skeleton. In section \ref{s3}, the
main tools among which the LDP for the Wiener process, the
continuity of the skeleton with respect to the potential on the
level sets of the rate function and exponential tail estimates are
presented. Section \ref{s4} is devoted to the proof of the uniform
LDP and Section \ref{s5} to an application to the study of
asymptotics for the blow-up time.

\section{Notations and preliminaries}\label{s2}
\subsection{Notations}
Throughout the paper the following notations will be used.\\
\indent The set of positive integers and positive real numbers are
denoted respectively by $\N^*$ and $\R_+^*$, while the set of real
numbers different from $0$ is denoted by $\R^*$.\vspace{0.3cm}

\indent For $p\in \N^*$, $\xLn^{p}(\R^d)$ is the classical Lebesgue
space of complex valued functions. For $k$ in $\N^*$,
$\xW^{k,p}(\R^d)$ is the associated Sobolev space of
$\xLn^{p}(\R^d)$ functions with partial derivatives up to level $k$,
in the sense of distributions, in $\xLn^{p}(\R^d)$. When $p=2$,
$\xHn^{s}(\R^d)$ denotes the fractional Sobolev space of tempered
distributions $v\in\mathcal{S}'(\R^d)$ such that the Fourier
transform $\hat{v}$ satisfies $(1+|\xi|^2)^{s/2}\hat{v}\in \L2$. It
is a Hilbert space. We denote by $\xHn^{s}(\R^d,\R)$ the space of
real-valued functions in $\xHn^{s}(\R^d)$. The space $\L2$ is
endowed with the inner product defined by
$(u,v)_{\L2}=\Re\int_{\R^d}u(x)\overline{v}(x)dx$.\vspace{0.3cm}

\indent If $I$ is an interval of $\R$, $(E,\|\cdot\|_E)$ a Banach
space and $r$ belongs to $[1,+\infty]$, then $\xLn^{r}(I;E)$ is the
space of strongly Lebesgue measurable functions $f$ from $I$ into
$E$ such that $t\rightarrow \|f(t)\|_E$ is in $\xLn^{r}(I)$. The
space $\xLn^{r}(\Omega;E)$, where
$\left(\Omega,\F,\mathbb{P}\right)$ is the probability space, is
defined similarly.\vspace{0.3cm}

\indent We recall that a pair $(r(p),p)$ of positive numbers is
called an admissible pair if $p$ satisfies $2\leq p<\frac{2d}{d-2}$
when $d>2$ ($2\leq p<+\infty$ when $d=2$ and $2\leq p\leq+\infty$
when $d=1$) and $r(p)$ is such that
$\frac{2}{r(p)}=d\left(\frac{1}{2}-\frac{1}{p}\right)$. For example
$(+\infty,2)$ is an admissible pair.\vspace{0.3cm}

Following the denomination of \cite{B1,B2}, the linear continuous
operator $\Phi$ is a $\gamma-$radonifying operator from a separable
Hilbert space $H$ into a Banach space $B$ if for any complete
orthonormal system $(e^H_j)_{j\in\N}$ of $H$ and
$(\gamma_j)_{j\in\N}$ a sequence of identically distributed and
independent Gaussian random variables on a probability space
$(\tilde{\Omega},\tilde{\mathcal{F}},\tilde{\mathbb{P}})$, we have
$$\|\Phi\|_{R(H,B)}^2=\tilde{\mathbb{E}}\left\|\sum_{j\in\N}\gamma_j\Phi
e^H_j\right\|_B^2<\infty.$$The space of such operators endowed with
the norm $\|\cdot\|_{R(H,B)}$ is a Banach space denoted by $R(H,B)$.
When $B$ is a Hilbert space $\tilde{H}$, it corresponds to the space
of Hilbert-Schmidt operators from $H$ into $\tilde{H}$. We denote by
$\mathcal{L}_2(H,\tilde{H})$ the space of Hilbert-Schmidt operators
from $H$ into $\tilde{H}$ endowed with the norm
$$\|\Phi\|_{\mathcal{L}_2(H,\tilde{H})}=\xtr\left(\Phi\Phi^*\right)
=\sum_{j\in\N}\|\Phi e^H_j\|_{\tilde{H}}^2,$$ where $\Phi^*$ denotes
the adjoint of $\Phi$ and $\xtr$ the trace. We denote by
$\mathcal{L}_2^{s,r}$ the corresponding space for
$H=\xHn^{s}(\R^d,\R)$ and $\tilde{H}=\xHn^{r}(\R^d,\R)$. The space
of linear continuous mappings from a Banach space $B$ into a Banach
space $\tilde{B}$ is denoted by
$\mathcal{L}_c\left(B,\tilde{B}\right)$. \vspace{0.3cm}

\indent When $A$ and $B$ are two Banach spaces, $A\cap B$, where the
norm of an element is defined as the maximum of the norm in $A$ and
in $B$, is a Banach space. Given an admissible pair $(r(p),p)$ and
$S$ and $T$ such that $0\leq S<T$, the space
$$X^{(S,T,p)}=\xCn\left([S,T];\xHone(\R^d)\right)\cap
\xL^{r(p)}\left(S,T;\xW^{1,p}(\R^d)\right),$$ denoted by $X^{(T,p)}$
when $S=0$ is of particular interest for the NLS
equation.\vspace{0.3cm}

We present hereafter the space of exploding paths, see \cite{EG} for
the main properties. We add a point $\Delta$ to the space $\H1$ and
embed the space with the topology such that its open sets are the
open sets of $\H1$ and the complement in $\H1\cup\{\Delta\}$ of the
closed bounded sets of $\H1$. This topology induces on $\H1$ the
topology of $\H1$. The set $\xCn([0,+\infty);\H1\cup\{\Delta\})$ is
then well defined. If $f$ belongs to
$\xCn([0,+\infty);\H1\cup\{\Delta\})$ we denote the blow-up time by
$$\mathcal{T}(f)=\inf\{t\in[0,+\infty):f(t)=\Delta\},$$
with the convention that $\inf\emptyset=+\infty$. We may now define
the set
$$\mathcal{E}(\H1)=\left\{f\in\xCn([0,+\infty);\H1\cup\{\Delta\}):
f(t_0)=\Delta\Rightarrow\forall t\geq t_0,\ f(t)=\Delta\right\},$$
it is endowed with the topology defined by the neighborhood basis
$$V_{T,\epsilon}(\varphi_1)=\left\{\varphi\in\mathcal{E}(\H1):
\mathcal{T}(\varphi)> T,\ \|\varphi_1-\varphi\|_{\xCn
\left([0,T];\H1\right)}\leq\epsilon\right\},$$ of $\varphi_1$ in
$\mathcal{E}(\H1)$ given $T<\mathcal{T}(\varphi_1)$ and
$\epsilon>0$.\vspace{0.3cm}

We define $\mathcal{A}(d)$ and $\tilde{\mathcal{A}}(d)$ as
$[2,+\infty)$ when $d=1$ or $d=2$ and respectively as
$\left[2,\frac{2(3d-1)}{3(d-1)}\right)$ and
$\left[2,\frac{2d}{d-1}\right)$ when $d\geq3$. The space
$\mathcal{E}_{\infty}$ is now defined for any $d$ in $\N^*$ by the
set of functions $f$ in $\mathcal{E}(\H1)$ such that for all
$p\in\mathcal{A}(d)$ and all $T\in [0,\mathcal{T}(f))$, $f$ belongs
to $\xL^{r(p)}\left(0,T;\xW^{1,p}(\R^d)\right)$. It is endowed with
the topology defined for $\varphi_1$ in $\mathcal{E}_{\infty}$ by
the neighborhood basis
$$W_{T,p,\epsilon}(\varphi_1)=\left\{\varphi\in\mathcal{E}_{\infty}
:\mathcal{T}(\varphi)>T,\
\|\varphi_1-\varphi\|_{X^{(T,p)}}\leq\epsilon \right\}$$ where
$T<\mathcal{T}(\varphi_1)$, $p\in\mathcal{A}(d)$ and $\epsilon>0$.
We may show using H\"older's inequality that it is a well defined
neighborhood basis. Also the space is a Hausdorff topological space
and thus we may consider applying generalizations of Varadhan's
contraction principle. If we denote again by
$\mathcal{T}:\mathcal{E}_{\infty}\rightarrow [0,+\infty]$ the
blow-up time, the mapping is measurable and lower
semicontinuous.\vspace{0.3cm}

We denote by $x\wedge y$ the minimum of the two real numbers $x$ and
$y$ and $x\vee y$ their maximum. We finally recall that a rate
function $I$ is a lower semicontinuous function and that a good rate
function $I$ is a rate function such that for every $c>0$,
$\left\{x:I(x)\leq c\right\}$ is a compact set.

\subsection{Properties of the group}\label{s22}
We recall in this section the main properties of the group that
are used in the article.\\
The group satisfies the decay estimates: for every $p\geq 2$, $t\neq
0$ and $u_0$ in $\xLn^{p'}(\R^d)$, where $p'$ is the conjugate
exponent of $p$, i.e. $\frac{1}{p}+\frac{1}{p'}=1$,
$$\|U(t)u_0\|_{\xW^{1,p}(\R^d)}\leq(4\pi|t|)^{-d\left(\frac{1}{2}
-\frac{1}{p}\right)}\|u_0\|_{\xW^{1,p'}(\R^d)}$$ and the following
Strichartz
estimates, see \cite{kato,SS}\\
\noindent $i/\ \ $For every $u_0$ in $\H1$ and $(r(p),p)$ an
admissible pair, the following linear mapping
\begin{align*}\xHone(\R^d)&\rightarrow \xCn(\R;\xHone(\R^d))\cap
\xLn^{r(p)}(\R;\xW^{1,p}(\R^d))\\
u_0&\mapsto (t\mapsto U(t)u_0).\end{align*} is continuous.\\
$ii/\ $For every $T$ positive and $(r(p),p)$ and $(r(q),q)$ two
admissible pairs, if $s$ and $\rho$ are the conjugate exponents of
$r(q)$ and $q$, the following linear mapping
\begin{align*}
\xLn^{s}(0,T;\xW^{1,\rho}(\R^d))&\rightarrow
\xCn([0,T];\xHone(\R^d))\cap
\xL^{r(p)}(0,T;\xW^{1,p}(\R^d))\\
f&\mapsto\Lambda f=\int_0^{\cdot}U(\cdot-s)f(s)ds\end{align*} is
continuous with a norm that does not depend on $T$.\vspace{0.3cm}

\subsection{Statistical properties of the noise}\label{s23}
We consider that $W$ is originated from a cylindrical Wiener process
on $\xLtwo(\R^d,\R)$, i.e. $W=\Phi W_c$ where $\Phi$ is a
$\gamma-$radonifying operator. It is known that for any orthonormal
basis $(e_j)_{j\in\N}$ of $\xLtwo(\R^d,\R)$ there exists a sequence
of real independent Brownian motions $(\beta_j)_{j\in\N}$ such that
the cylindrical Wiener process may be written as the following
expansion on a complete orthonormal system
$\left(e_j\right)_{j\in\N}$ of $\xLtwo(\R^d,\R)$
$$W_c(t,x,\omega)=\sum_{j\in\N}\beta_j(t,\omega) e_j(x).$$
Its direct image in every Banach space $B$ such that the injection
of $\xLtwo(\R^d,\R)$ into $B$ is a $\gamma-$radonifying injection is
a {\it bona fide} Wiener process. It is also possible to define the
Wiener process as originated from a cylindrical Wiener process on
the reproducing kernel Hilbert space (RKHS) of the measure $\mu$,
the law of $W(1)$, in such a case the mapping $\Phi$ is an
injection. The assumption that $\Phi$ is $\gamma-$radonifying is
exactly the assumption needed to transfer by direct image a
cylindrical centered Gaussian measure on a separable Hilbert space
to a {\it bona fide} $\sigma-$additive measure $\mu$ on the image
Banach space. Remark also that this assumption is too severe to
allow the noise to be homogeneous in space. In \cite{dBD1}, the
assumption on the noise is that $\Phi\in
R\left(\xLtwo(\R^d,\R),\xW^{1,\alpha}(\R^d)\right)\cap\mathcal{L}_2^{0,1}$,
where $\alpha>2d$. In particular the noise is a real noise. Stronger
assumptions on the spacial correlations of the Wiener process than
in the additive case, see \cite{dBD1,EG}, are imposed in order that
the stochastic convolution of the product make sense in the space
considered for the fixed point.\vspace{0.3cm}

In this article we impose the extra assumption (A) \vspace{0.2cm}

\begin{center}
for some $s>\frac{d}{4}+1$, $\Phi\in \mathcal{L}_2^{0,s}$.
\end{center}
\vspace{0.2cm}

It is used twice in the proof of the LDP. It is used first to prove
the continuity with respect to the potential on the level sets of
the rate function of the sample path LDP for the Wiener process. It
is needed to produce a process in a real separable Banach space
embedded in $\xW^{1,\infty}(\R^d)$. Note that in that respect it
would be enough to impose that $\Phi\in R(\L2,
\xHone(\R^d,\R)\cap\xW^{1+s,\beta}(\R^d))$ for any
$s\beta>\frac{d}{2}$. It is used a second time in the proof of the
exponential tail estimates. In that case we also need that
$\xHn^{s}(\R^d,\R)$ is a Hilbert space.\vspace{0.3cm}

Since the operator $\Phi$ belongs to $\mathcal{L}_2^{0,s}$ and thus
to $\mathcal{L}_2^{0,0}$ it may be defined through a kernel. This
means that for any square integrable function $u$,
$$\Phi u(x) = \int_{\R^d}
\mathcal{K}(x,y) u(y) dy.$$ Now for $(t,s)\in\R^+$,
$(x,z)\in(\R^d)^2$, the correlation of our real Gaussian noise
$\frac{\partial}{\partial t}W$ is formally given by
$$\E\left[\frac{\p W}{\p t}(t+s,x+z)\frac{\p W}{\p
t}(t,x)\right]=\delta_0(s)\otimes c(x,z),$$ where the distribution
$c(x,z)$ is indeed the
 $\L2$ function
$$c(x,z)=\int_{\R^d}\mathcal{K}(x+z,u)\mathcal{K}(x,u)du.$$

In the following we assume that the probability space is endowed
with the filtration $\mathcal{F}_t=\mathcal{N}\cup\s\{W_s,0\leq
s\leq t\}$ where $\mathcal{N}$ denotes the $\mathbb{P}-$null sets.

\subsection{The law of the solution in the space of exploding paths}
In the following we restrict ourselves to the case where
$$\left\{\begin{array}{cl}
&\frac{1}{2}\leq\sigma\hspace{1.3cm} \mbox{if}\ d=1,2,\\
&\frac{1}{2}\leq\sigma<\frac{2}{d-2}\ \ \mbox{if}\ d\geq3.
\end{array}\right.$$
\begin{prpstn}\label{prpstn1}
The solution $u^{u_0}$ defines a random variable with values in
$\mathcal{E}_{\infty}$.
\end{prpstn}
\begin{proof}
The proof of existence and uniqueness of a maximal solution in
\cite{dBD1}, in the case where $\sigma\geq\frac{1}{2}$, is obtained
as follows. Take $p\in\mathcal{A}(d)$ and $T>0$. The mapping
\begin{align*}
\mathcal{F}_{R,u_0}(u)(t)=&U(t)u_0-i\lambda\int_0^tU(t-s)
\left(\theta\left(\frac{\|u\|_{X^{(s,p)}}}{R}\right)|u(s)|
^{2\sigma}u(s)\right)ds\\
&-i\int_0^tU(t-s)\left(u(s)dW(s)\right)-\frac{1}{2}\int_0^t
U(t-s)\left(u(s)F_{\Phi}\right)ds
\end{align*}
is a contraction in $\xLn^{r(p)}\left(\Omega,X^{(T^*,p)}\right)$
provided $T^*$, depending on $p$, $\|u_0\|_{\H1}$ and $R$, is small
enough. The fixed point is denoted by $u_R^{u_0}$. The solution can
be extended to the whole interval $[0,T]$ and is a measurable
mapping from $(\Omega,\mathcal{F})$ to $X^{(T,p)}$ with its Borel
$\sigma-$field. The blow-up time is defined as the increasing limit
of the approximate blow-up time, see also \cite{EG},
$$\tau_R=\inf\{t\in\R^+:\|u_R^{u_0}\|_{X^{(t,p)}}\geq R\}.$$
The solution $u^{u_0}$ is then such that $u^{u_0}=u_R^{u_0}$ on
$[0,\tau_R]$. Finally it is proved that the limit of the approximate
blow-up time corresponds indeed to the blow-up of the $\xHone(\R^d)$
norm. Consequently, the solution is an element of
$\mathcal{E}_{\infty}$. The topology of $\mathcal{E}_{\infty}$ is
defined by a countable basis of neighborhoods, see \cite{EG}, thus
it suffices to show that for $\varphi_1\in \mathcal{E}_{\infty}$,
$T<\mathcal{T}(\varphi_1)$, $p\in\mathcal{A}(d)$ and $\epsilon>0$,
$\left\{u^{u_0}\in W_{T,p,\epsilon}(\varphi_1)\right\}$ is an
element of $\mathcal{F}$. Since $T<\mathcal{T}(\varphi_1)$ there
exists $R>0$ such that $\|\varphi_1\|_{X^{(T,p)}}\leq R$ and every
$\varphi$ in $W_{T,p,\epsilon}(\varphi_1)$ satisfies
$\|\varphi\|_{X^{(T,p)}}\leq R+\epsilon$. The conclusion follows
from the fact that
$$\left\{u^{u_0}\in W_{T,p,\epsilon}(\varphi_1)\right\}=
\left\{u_{R+\epsilon}^{u_0}\in
W_{T,p,\epsilon}(\varphi_1)\right\}=\left\{\|\varphi_1-
u_{R+\epsilon}^{u_0}\|_{X^{(T,p)}} \leq \epsilon\right\}.$$
\end{proof}
We denote by $\mu^{u^{u_0}}$ its law and, for $\epsilon>0$, by
$\mu^{u^{\epsilon,u_0}}$ the law in $\mathcal{E}_{\infty}$ of the
mild solution of $$\left\{\begin{array}{lc}
  i \xdif u^{\epsilon,u_0}=\left(\Delta u^{\epsilon,u_0} + \lambda
|u^{\epsilon,u_0}|^{2\s}u^{\epsilon,u_0} -
\frac{i\epsilon}{2}F_{\Phi}u^{\epsilon,u_0}\right)\xdif t +
\sqrt{\epsilon}u^{\epsilon,u_0} \xdif W.\\
u^{\epsilon,u_0}(0)=u_0
\end{array}\right.$$

\subsection{The uniform large deviation principle}
The rate function of the LDP involves the deterministic control
problem
$$\left\{\begin{array}{cl}
&i\frac{\xdif}{\xdif t}u=\Delta u +\lambda |u|^{2\sigma}u + u\Phi h,\\
&u(0)=u_0\ \in \xHone(\R^d),\\
&h\in \xLtwo\left(0,+\infty;\xLtwo(\R^d)\right).
\end{array}\right.$$
We may write the mild solution, also called skeleton, as
$$S^c(u_0,h)=U(t)u_0-i\int_0^tU(t-s)\left[S^c(u_0,h)(s)\left(
\lambda|S^c(u_0,h)(s)|^{2\sigma}+\Phi h(s)\right)\right]ds.$$ If we
replace $\Phi h$ by $\frac{\p f}{\p t}$ where $f$ belongs to
$\xHn^1_0\left([0,+\infty);\xHn^{s}(\R^d,\R)\right)$ which is the
subspace of $\xCn\left([0,+\infty);\xHn^{s}(\R^d,\R)\right)$ of
functions null at time $0$, square integrable in time and with
square integrable in time time derivative, we write
$$S(u_0,f)=U(t)u_0-i\int_0^tU(t-s)\left[S(u_0,f)(s)\left(\lambda|S(u_0,f)(s)|
^{2\sigma} +\frac{\p f}{\p s}\right)\right]ds.$$ In the following we
denote by $K\subset\subset \xHone(\R^d)$ the fact that $K$ is a
compact set of $\xHone(\R^d)$ and by $Int(A)$ the interior of $A$.
\begin{thrm}
The family of measures
$\left(\mu^{u^{\epsilon,u_0}}\right)_{\epsilon>0}$ satisfies a
uniform LDP on $\mathcal{E}_{\infty}$ of speed $\epsilon$ and good
rate function \begin{align*}
I^{u_0}(w)&=\inf_{f\in\xCn\left([0,+\infty);\xHn^{s}(\R^d,\R)\right):w=S(u_0,f)}
I^W(f)\\
&=\frac{1}{2}\inf_{h\in\xLtwo(0,+\infty;\xLtwo(\R^d)):w=S^c(u_0,h)}\|h\|
_{\xLtwo\left(0,+\infty;\xLtwo(\R^d)\right)}^2,
\end{align*} where $I^W$ is the rate function of the sample path LDP for
the Wiener process, i.e. $\forall K\subset\subset\xHone(\R^d),\
\forall A\in\mathcal{B}\left(\mathcal{E}_{\infty}\right),$
$$-\sup_{u_0\in K}\inf_{w\in
Int(A)}I^{u_0}(w)\leq\underline{\lim}_{\epsilon\rightarrow0}\epsilon\log
\inf_{u_0\in K}\mathbb{P}\left(u^{\epsilon,u_0}\in A\right)$$
$$\leq\overline{\lim}_{\epsilon\rightarrow0}\epsilon\log\sup_{u_0\in
K}\mathbb{P}\left(u^{\epsilon,u_0}\in
A\right)\leq-\inf_{w\in\overline{A},u_0\in K}I^{u_0}(w).$$
\end{thrm}
\begin{rmrk}
Writing $\frac{\partial h}{\partial t}$ instead of $h$ in the
optimal control problem leads to a rate function consisting in the
minimisation of $\frac{1}{2}\|h\|^2_{H_0^1(0,+\infty;\L2)}$.
Specifying only the law $\mu$ of $W(1)$ and dropping $\Phi$ in the
control problem would lead to a rate function consisting in the
minimisation of $\frac{1}{2}\|h\|^2_{H_0^1(0,+\infty;H_{\mu})}$,
where $H_{\mu}$ stands for the RKHS of $\mu$.
\end{rmrk}

\section{The main tools}\label{s3}
In the two next sections $C$ is a constant which may have a
different value from line to line and also within the same line.
\subsection{Large deviations for the Wiener process}
\begin{prpstn}\label{prpstn2}
The family of laws of $\left(\sqrt{\epsilon}W\right)_{\epsilon>0}$
satisfies in $\xCn([0,+\infty);\xHn^{s}(\R^d,\R))$ a LDP of speed
$\epsilon$ and good rate function
$$I^{W}(f)=\frac{1}{2}\inf_{h\in\xLtwo\left(0,+\infty;\xLtwo(\R^d)\right):
f=\int_0^{\cdot}\Phi
h(s)ds}\|h\|_{\xLtwo\left(0,+\infty;\xLtwo(\R^d) \right)}^2.$$
\end{prpstn}
\begin{proof}
The result follows from the general LDP for Gaussian measures on a
real separable Banach space, see \cite{DS1}, the fact that the laws
of the restrictions of the paths $\sqrt{\epsilon}W$ on
$\xCn\left([0,T];\xHn^s(\R^d,\R)\right)$ and
$\xLtwo\left(0,T;\xHn^s(\R^d,\R)\right)$ have same RKHS and
Dawson-G\"artner's theorem which allows to deduce the LDP on
$\xCn\left([0,+\infty);\xHn^s(\R^d,\R)\right)$ with the topology of
the uniform convergence on the compact sets of $[0,+\infty)$.
\end{proof}
In the following we denote by $C_a$ the set
\begin{center}
\begin{tabular}{rl}
$C_a$&$=\left\{f\in\xCn([0,+\infty);\xHn^{s}(\R^d,\R)):I^W(f)\leq
a\right\},$\\
&$=\left\{f\in\xCn([0,+\infty);\xHn^{s}(\R^d,\R)): f(0)=0,\ \frac{\p
f}{\p t}\in\xim\Phi\right.$\\
&$\ \ \ \left.\mbox{and}\
\left\|\Phi^{-1}_{|(\xker\Phi)^{\perp}}\frac{\p f}{\p
t}\right\|_{\xLtwo(0,+\infty;\L2)}\leq \sqrt{2a}\right\}.$
\end{tabular}
\end{center}

\subsection{Continuity of the skeleton with respect to the initial
data and control on the level sets of the rate function of the
Wiener process} The continuity with respect to the control on the
level sets of the rate function of the Wiener process is used herein
to prove that the rate function is a good rate function and to prove
the lower bound of the LDP, see \ref{s41}. In that respect our proof
is closer to the proof of \cite{CW,CM}. The authors of \cite{CR,P1}
use some slightly different arguments and do not prove this
continuity but prove separately the compactness of the level sets of
the rate function, the lower and upper bounds of the LDP using the
usual characterization in metric spaces. Remark also that in the
applications of the LDP, see for example Section \ref{s5}, we need
to compute infima of the rate function, or of the rate function of a
LDP deduced by contraction, on particular sets and the continuity
proves to be very useful. We prove the stronger continuity with
respect to the control and initial datum as suggested in \cite{CM}
but we will not need the continuity with respect to the initial
datum in the proof of the uniform LDP.

\begin{prpstn}\label{prpstn3}
For every $u_0\in\H1$, $a>0$ and $f\in C_a$, $S(u_0,f)$ exists and
is uniquely defined. Moreover, it is a continuous mapping from
$\xHone(\R^d)\times C_a$ into $\mathcal{E}_{\infty}$, where $C_a$
has the topology induced by that of
$\xCn([0,+\infty);\xHn^{s}(\R^d,\R))$.
\end{prpstn}
\begin{proof}
Let $\mathfrak{F}$ denote the mapping such that
$$\mathfrak{F}(u,u_0,f)=U(t)u_0-i\int_0^tU(t-s)\left[u(s)
\left(\lambda|u(s)|^{2\sigma}+\frac{\p f}{\p s}\right)\right]ds.$$
Let $a$ and $r$ be positive, $f$ be in $C_a$ and $u_0$ be such that
$\|u_0\|_{\H1}\leq r$, set $R=2cr$ where $c$ is the norm of the
linear continuous mapping of the $i/$ of the Strichartz estimates.
From $i/$ and $ii/$ of the Strichartz estimates along with
H\"{o}lder's inequality, the Sobolev injections and the continuity
of $\Phi$, for any $T$ positive, $p$ in $\mathcal{A}(d)$ and $u$ and
$v$ in $X^{(T,p)}$ and any $\nu$ in
$\left(0,1-\frac{\sigma(d-2)}{2}\right)$, which is a well defined
interval since $\sigma<\frac{2}{d-2}$, there exists $C$ positive
such that \begin{center}
\begin{tabular}{l}$\|\mathfrak{F}(u,u_0,f)\|_{X^{(T,p)}}$\\
$\leq c\|u_0\|_{\H1}+CT^{\nu}\|u\|_{X^{(T,p)}}^{2\sigma+1}
+CT^{\frac{1}{2}-\frac{1}{r(p)}}\|u\|_{\xCn\left([0,T];\H1\right)}
\left\|\frac{\p
f}{\p s}\right\|_{\xLtwo\left(0,T;\xW^{1,\frac{r(p)d}{2}}(\R^d)\right)}$\\
$\leq c\|u_0\|_{\H1}+CT^{\nu}\|u\|_{X^{(T,p)}}^{2\sigma+1}+C\sqrt{a}
T^{\frac{1}{2}-\frac{1}{r(p)}}\|u\|_{X^{(T,p)}}$.
\end{tabular}\end{center}
Similarly we obtain
\begin{center}
\begin{tabular}{l}$\|\mathfrak{F}(u,u_0,f)-\mathfrak{F}(v,u_0,f)\|_{X^{(T,p)}}$\\
$\leq
C\left[T^{\nu}\left(\|u\|_{X^{(T,p)}}^{2\sigma}+\|v\|_{X^{(T,p)}}
^{2\sigma}\right)+T^{\frac{1}{2}-\frac{1}{r(p)}}\sqrt{a}\right]
\|u-v\|_{X^{(T,p)}}.$\end{tabular}\end{center} Thus, for
$T=T^*_{r,a,p}$ small enough depending on $r$, $a$ and $p$, the ball
centered at $0$ of radius $R$ is invariant and the mapping
$\mathfrak{F}(\cdot,u_0,f)$ is a $\frac{3}{4}-$contraction. We
denote by $S^0(u_0,f)$ the unique fixed point of
$\mathfrak{F}(\cdot,u_0,f)$ in
$X^{(T^*_{r,a,p},p)}$.\\
\indent Also when $T$ is positive, we can solve the fixed point
problem on any interval $\left[kT^*_{r,a,p},(k+1)T^*_{r,a,p}\right]$
with $1\leq k\leq \left\lfloor\frac{T}{T^*_{r,a,p}}\right\rfloor$.
The fixed point is denoted by $S^k(u_k,f)$ where
$u_k=S^{k-1}(u_{k-1},f)(kT^*_{r,a,p})$, as long as
$\|u_k\|_{\H1}\leq r$. Existence and uniqueness of a maximal
solution $S(u_0,f)$ follows. It coincides with $S^k(u_k,f)$ on the
above intervals when it is defined. We may also show that the
blow-up time corresponds to the blow-up of the $\H1-$norm, thus
$S(u_0,f)$ is an element
of $\mathcal{E}_{\infty}$.\\
\indent We shall now prove the continuity. Take $u_0$ in $\H1$, $a$
positive and $f$ in $C_a$. From the definition of the neighborhood
basis of the topology of $\mathcal{E}_{\infty}$, it is enough to see
that for $\epsilon$ positive, $T<\mathcal{T}\left(S(u_0,f)\right)$,
and $p\in\mathcal{A}(d)$, there exists $\eta$ positive such that for
every $\tilde{u}_0$ in $\H1$ and $g$ in $C_a$ satisfying
$\|u_0-\tilde{u}_0\|_{\H1}+\|f-g\|_{\xCn\left([0,T];\xHn^s(\R^d,\R)\right)}\leq
\eta$ then
$\left\|S(u_0,f)-S(\tilde{u}_0,g)\right\|_{X^{(T,p)}}\leq\epsilon$.
We set
$$r=\left\|S(u_0,f)\right\|_{X^{(T,p)}}+1,\
N=\left\lfloor\frac{T}{T^*_{r,a,p}}\right\rfloor,\
\delta_{N+1}=\frac{\epsilon}{N+1}\wedge 1,$$ and define for
$k\in\left\{0,...,N\right\},\ \delta_k$ and $\eta_k$ such that
$0<\delta_k<\delta_{k+1}$, $0<\eta_k<\eta_{k+1}<1$ and
$$\left\|S^{k+1}(u_k,f)-S^{k+1}(\tilde{u}_k,g)\right\|
_{X^{(kT^*_{r,a,p},(k+1)T^*_{r,a,p},p)}}\leq\delta_{k+1},$$ if
$$\|u_k-\tilde{u}_k\|_{\H1}+\|f-g\|_{\xCn\left([0,\infty);\xHn^s(\R^d,\R)
\right)}\leq\eta_{k+1}.$$ It is possible to choose
$\left(\delta_k\right)_{k=0,...,N}$ as long as we prove the
continuity for $k=0$. The maximal solution $S(\tilde{u}_0,g)$ then
necessarily satisfies $\mathcal{T}\left(S(\tilde{u}_0,g)\right)>T$.
We conclude setting $\eta=\eta_1$ and
using the triangle inequality.\\
\indent We now prove the continuity for $k=0$. First note that
$\|u_0\|_{\H1}\leq r$ and $\|\tilde{u}_0\|_{\H1}\leq r$ since
$\eta<1$, as a consequence if $\Upsilon$ denotes\\
$\left\|S^0(u_0,f)-S^0(\tilde{u}_0,g)\right\|_{X^{\left(T^*_{r,a,p},p\right)}}$
 we obtain that
\begin{align*}
\Upsilon=&\left\| \mathfrak{F}(S^0(u_0,f),u_0,f)
-\mathfrak{F}(S^0(\tilde{u}_0,g),\tilde{u}_0,g)\right\|_{X^{\left(T^*_{r,a,p},p\right)}}\\
&\leq\left\|\mathfrak{F}(S^0(u_0,f),u_0,f)-\mathfrak{F}(S^0(u_0,f),\tilde{u}_0,g)
\right\|_{X^{\left(T^*_{r,a,p},p\right)}}\\
&\ \ \
+\left\|\mathfrak{F}(S^0(u_0,f),\tilde{u}_0,g)-\mathfrak{F}(S^0(\tilde{u}_0,g),
\tilde{u}_0,g) \right\|_{X^{\left(T^*_{r,a,p},p\right)}},
\end{align*}
thus since $\mathfrak{F}(\cdot,\tilde{u}_0,g)$ is a
$\frac{3}{4}-$contraction, \begin{align*}\Upsilon&\leq
4\left\|\mathfrak{F}(S^0(u_0,f),u_0,f)-\mathfrak{F}(S^0(u_0,f),\tilde{u}_0,g)\right\|
_{X^{\left(T^*_{r,a,p},p\right)}}\\
&\leq
4\left(c\|u_0-\tilde{u}_0\|_{\H1}+\left\|\int_0^{\cdot}U(\cdot-s)S^0(u_0,f)\frac{\p
(f-g)}{\p
s}(s)ds\right\|_{X^{\left(T^*_{r,a,p},p\right)}}\right).\end{align*}
Using H\"older's inequality and taking $p<\tilde{p}$ we can bound
the second term, denoted by $T_2$, by
\begin{align*}
T_2\leq&\left(\left\|\int_0^{\cdot}U(\cdot-s)S^0(u_0,f)\frac{\p
(f-g)}{\p
s}(s)ds\right\|_{\xCn\left(\left[0,T^*_{r,a,p}\right];\H1\right)}^{\theta}\right.\\
&\left.\times\left\|\int_0^{\cdot}U(\cdot-s)S^0(u_0,f)\frac{\p
(f-g)}{\p
s}(s)ds\right\|_{\xLn^{r(\tilde{p})}\left(0,T^*_{r,a,p};\xW^{1,\tilde{p}}
(\R^d)\right)}^{1-\theta}\right)\\
&\vee\left\|\int_0^{\cdot}U(\cdot-s)S^0(u_0,f)\frac{\p (f-g)}{\p
s}(s)ds\right\|_{\xCn\left(\left[0,T^*_{r,a,p}\right];\H1\right)},
\end{align*}
where $\theta=\frac{\tilde{p}-p}{\tilde{p}-2}$. From the $ii/$ of
the Strichartz estimates we can bound the second term of the product
by
$C\left(\sqrt{a}\left(T^*_{r,a,p}\right)^{\frac{1}{2}-\frac{1}{r(\tilde{p})}}R
\right)^{1-\theta}$. It is now enough to show that when $g$ is close
enough to $f$,
the first term in the above can be made arbitrarily small.\\
\indent Take $n\in \mathbb{N}$ and set for
$i\in\left\{0,...,n\right\},\ t_i=\frac{iT^*_{r,a,p}}{n}$ and
$S^{0,n}(u_0,f)=U(t-t_i)\left(S(u_0,f)(t_i)\right)$ when $t_i\leq
t<t_{i+1}$. As it has been done previously we can compute the
following bound,
\begin{center}
\begin{tabular}{l}
$\left\|\int_0^{\cdot}U(\cdot-s)\left(S^0(u_0,f)-S^{0,n}(u_0,f)\right)\frac{\p
(f-g)}{\p
s}(s)ds\right\|_{\xCn\left(\left[0,T^*_{r,a,p}\right];\H1\right)}$\\
$\leq
C\sqrt{a}\left(T^*_{r,a,p}\right)^{\frac{1}{2}-\frac{1}{r(p)}}\left\|S^0(u_0,f)-S^{0,n}(u_0,f)
\right\|_{\xCn\left(\left[0,T^*_{r,a,p}\right];\H1\right)}$\\
$\leq
C\sqrt{a}\left(T^*_{r,a,p}\right)^{\frac{1}{2}-\frac{1}{r(p)}}$\\
$\sup_{i\in\{0,...,n-1\}}
\left\|\int_{t_i}^{\cdot}U(\cdot-s)\left[S^0(u_0,f)\left(\lambda
\left|S^0(u_0,f)\right|^{2\sigma}+\frac{\p f}{\p
s}\right)\right](s)ds\right\|_{\xCn\left([t_i,t_{i+1}];\H1\right)}$\\
$\leq
C\sqrt{a}\left(T^*_{r,a,p}\right)^{\frac{1}{2}-\frac{1}{r(p)}}\left[R^{2\sigma+1}
\left(\frac{T^*_{r,a,p}}{n}\right)^{\nu}+\left(\frac{T^*_{r,a,p}}{n}
\right)^{\frac{1}{2}-\frac{1}{r(p)}}R\sqrt{a}\right],$
\end{tabular}\end{center}
which can be made arbitrarily small for large $n$. Finally it
remains to bound the
$\xCn\left(\left[0,T^*_{r,a,p}\right];\H1\right)-$norm of
\begin{center}
\begin{tabular}{l} $\int_0^tU(t-s)S^{0,n}(u_0,f)\frac{\p
(f-g)}{\p
s}(s)ds$\\
$=\sum_{i=0}^{n-1}\int_{t_i\wedge t}^{t_{i+1}\wedge t}U(t-t_i\wedge
t)\left(S^0(u_0,f)(t_i\wedge t)\frac{\p (f-g)}{\p
s}(s)\right)ds$\\
$=\sum_{i=0}^{n-1}U(t-t_i\wedge t)\left[S^0(u_0,f)(t_i\wedge
t)\left((f-g)(t_{i+1}\wedge t)-(f-g)(t_i\wedge t)\right)\right]$
\end{tabular}\end{center}
since, from the Sobolev injections and the fact that $U(t-t_i\wedge
t)$ is a group on $\H1$, for any $i\in\left\{0,...,n-1\right\}$ and
$v$ in $\xHn^s(\R^d,\R)$,
$$\left\|U(t-t_i\wedge t)\left(S^0(u_0,f)(t_i\wedge
t)v\right)\right\|_{\H1}\leq C \left\|S^0(u_0,f)(t_i\wedge
t)\right\|_{\H1}\|v\|_{\xHn^s(\R^d,\R)},$$ and from the fact that
$\frac{\p (f-g)}{\p s}$ is Bochner integrable. We finally obtain
that an upper bound can be written as
$rCn\|f-g\|_{\xCn\left([0,T];\xHn^s(\R^d,\R)\right)}$, which for
fixed $n$ can be made arbitrarily small taking $g$ sufficiently
close to $f$.
\end{proof}

\subsection{Exponential tail estimates}
In the following we denote by $c\left(\frac{r(p)d}{2}\right)$ and
$c(\infty)$ the norms of the linear continuous embeddings
$\xHn^s(\R^d,\R)\subset\xW^{1,\frac{r(p)d}{2}}(\R^d,\R)$ and
$\xHn^s(\R^d,\R)\subset\xW^{1,\infty}(\R^d,\R)$. Exponential tail
estimates for the $\xLinfty-$norm of stochastic integrals in
$\xLn^p(0,1)$ is proved in \cite{CW} but here we need the following
\begin{lmm}\label{l1}Assume that $\xi$ is a point-predictable process,
that $p$ is such that $p\in\tilde{\mathcal{A}}(d)$, that $T>0$ and
that there exists $\eta>0$ such that
$\|\xi\|_{\xCn\left([0,T];\xHone(\R^d)\right)}^2\leq\eta$ a.s., then
for every $t\in[0,T]$, and $\delta>0$,
$$\mathbb{P}\left(\sup_{t_0\in[0,T]}\left\|\int_0^{t_0}U(t-s)\xi(s)dW(s)
\right\|_{\xW^{1,p}(\R^d)}\geq \delta\right)\leq
\exp\left(1-\frac{\delta^2} {\kappa(\eta)}\right),$$ where
$$\kappa(\eta)=\frac{4c\left(\frac{r(p)d}{2}\right)^2T^{1-\frac{4}{r(p)}}(d+1)
(d+p)\|\Phi\|_{\mathcal{L}_2^{0,s}}^2}{1-\frac{4}{r(p)}}\eta.$$
\end{lmm}
\begin{proof}
Let us denote by
$g_a(f)=\left(1+a\|f\|_{\xW^{1,p}(\R^d)}^p\right)^{\frac{1}{p}}$ the
real-valued function parametrized by $a>0$ and by $M$ the martingale
defined by $M(t_0)=\int_0^{t_0}U(t-s)\xi(s)dW(s)$. The function
$g_a$ is twice Fr\'echet differentiable, the first and second
derivatives at point $M(t)$ in the direction $h$ are denoted by
$Dg_a(M(t)).h$ and $D^2g_a(M(t),M(t)).(h,h)$, they are continuous.
Also, the second derivative is uniformly continuous on the bounded
sets. By the It\^o formula the following decomposition holds
$$g_a(M(t_0))=1+E_a(t_0)+R_a(t_0),$$
where $E_a(t_0)$ is equal to
$$\int_0^{t_0}Dg_a(M(s)).U(t-s)\xi(s)dW(s)-\frac{1}{2}
\int_0^t\|Dg_a(M(s)).U(t-s)\xi(s)\|_{\mathcal{L}_2(\L2,\R)}^2ds$$
and $R_a(t_0)$ to
\begin{align*}
\frac{1}{2}&\left(\int_0^{t_0}\|Dg_a(M(s)).U(t-s)\xi(s)\|
_{\mathcal{L}_2(\L2,\R)}^2ds\right.\\
&\left.+\int_0^{t_0}\sum_{j\in\N}D^2g_a(M(s),M(s)).(U(t-s)\xi(s)\Phi
e_j,U(t-s)\xi(s)\Phi e_j)ds\right),
\end{align*}
where $(e_j)_{j\in\N}$ is a complete orthonormal system of $\L2$. We
denote by
$$q(u,h)=\Re\left[\int_{\R^d}\overline{u}|u|^{p-2}hdx+\sum_{k=1}^d
\int_{\R^d}\overline{\p_{x_j}u}|\p_{x_j}u|^{p-2}\p_{x_j}hdx\right],$$
then
$$Dg_a(u).h=a\left(1+a\|u\|_{\xW^{1,p}(\R^d)}^p\right)^{\frac{1}{p}-1}q(u,h)$$
and
\begin{align*}
&D^2g_a(u,u).(h,g)=a^2(1-p)\left(1+a\|u\|_{\xW^{1,p}(\R^d)}^p
\right)^{\frac{1}{p}-2}q(u,h)q(u,g)\\
&+a\left(1+a\|u\|_{\xW^{1,p}(\R^d)}^p\right)^{\frac{1}{p}-1}
\left[\frac{p}{2}\Re\int_{\R^d}\left(|u|^{p-2}\overline{g}h
+\sum_{k=1}^d|\p_{x_j}u|^{p-2}\overline{\p_{x_j}g}\p_{x_j}h\right)dx\right.\\
&\left.+\frac{p-2}{2}\int_{\R^d}\left(\overline{u}^2|u|^{p-4}gh
+\overline{\p_{x_j}u}^2|\p_{x_j}u|^{p-4}\p_{x_j}g\p_{x_j}h\right)dx\right].
\end{align*}
From the series expansion of the Hilbert-Schmidt norm along with
H\"older's inequality, we obtain that $R_a(t_0)$ is bounded above by
\begin{align*}
&\frac{a(d+1)}{2}\int_0^{t_0}\left[a(d+1)\sum_{j\in\N}
\|U(t-s)\xi(s)\Phi e_j\|_{\xW^{1,p}(\R^d)}^2\left(\frac{\|M(s)\|
_{\xW^{1,p}(\R^d)}}{\left(1+a\|M(s)\|_{\xW^{1,p}(\R^d)}^p
\right)^{\frac{1}{p}}}\right)^{2(p-1)}\right.\\
&\left.-(p-1)a(d+1)\left(1+a\|M(s)\|_{\xW^{1,p}(\R^d)}^p
\right)^{\frac{1}{p}-2}\sum_{j\in\N}q\left(M(s),U(t-s)\xi(s)\Phi
e_j\right)^2\right.\\
&\left.(p-1)\|M(s)\|_{\xW^{1,p}(\R^d)}^{p-2}
\left(1+a\|M(s)\|_{\xW^{1,p}(\R^d)}^p\right)^{\frac{1}{p}-1}
\sum_{j\in\N}\|U(t-s)\xi(s)\Phi e_j\|_{\xW^{1,p}(\R^d)}^2\right]ds.
\end{align*}
Since the term in parenthesis in the first part is a decreasing
function of $\|M(s)\|_{\xW^{1,p}(\R^d)}$, the second term is non
positive and the following calculation holds
\begin{align*}
&\|M(s)\|_{\xW^{1,p}(\R^d)}^{p-2}\left(1+a\|M(s)\|_{\xW^{1,p}(\R^d)}^p
\right)^{\frac{1}{p}-1}\\
&= a^{\frac{2}{p}-1}\left(a\|M(s)\|_{\xW^{1,p}(\R^d)}^p
\right)^{1-\frac{2}{p}}\left(1+a\|M(s)\|_{\xW^{1,p}(\R^d)}^p\right)^{\frac{1}{p}-1}\\
&\leq a^{\frac{2}{p}-1}\left(1+a\|M(s)\|_{\xW^{1,p}(\R^d)}^p
\right)^{1-\frac{2}{p}}\left(1+a\|M(s)\|_{\xW^{1,p}(\R^d)}^p\right)^{\frac{1}{p}-1}\leq
a^{\frac{2}{p}-1},
\end{align*}
we obtain that
$$R_a(t_0)\leq\frac{(d+1)(d+p)a^{\frac{2}{p}}}{2}\int_0^{t_0}
\sum_{j\in\N}\|U(t-s)\xi(s)\Phi e_j\|_{\xW^{1,p}(\R^d)}^2ds.$$
Finally from the decay estimates, see Section \ref{s22}, H\"older's
inequality and the Sobolev injections we obtain that for any
$t_0\in[0,T]$,
$$R_a(t_0)\leq\frac{2(d+1)(d+p)a^{\frac{2}{p}}}{4}c
\left(\frac{r(p)d}{2}\right)^2\|\Phi\|_{\mathcal{L}_2^{0,s}}^2
\eta\int_0^T|t-s|^{-\frac{4}{r(p)}}ds,$$ the integral is finite
since we have made the assumption that $p<\frac{2d}{d-1}$ and thus
$$R_a(t_0)\leq\frac{\kappa(\eta)a^{\frac{2}{p}}}{4}.$$
Also calculation using the fact that
$\left(\exp\left(E_a(t_0)\right)\right)_{t_0\in[0,T]}$ is a
martingale, indeed the Novikov condition is satisfied from the above
and Doob's inequality leads
\begin{align*}
&
\mathbb{P}\left(\sup_{t_0\in[0,T]}\left\|\int_0^{t_0}U(t-s)\xi(s)dW(s)
\right\|_{\xW^{1,p}(\R^d)}\geq \delta\right)\\
&=\mathbb{P}\left(\sup_{t_0\in[0,T]}\exp\left(
g_a(M(t_0))\right)\geq
\exp\left((1+a\delta^p)^{\frac{1}{p}}\right)\right)\\
&\leq\mathbb{P}\left(\sup_{t_0\in[0,T]}\exp
\left(E_a(t_0)\right)\geq\exp\left((1+a\delta^p)^{\frac{1}{p}}-1
-\frac{\kappa(\eta)a^{\frac{2}{p}}}{4}\right)\right)\\
&\leq\exp\left(-(1+a\delta^p)^{\frac{1}{p}}+1+\frac{\kappa(\eta)
a^{\frac{2}{p}}}{4}\right)\\
&\leq
e\exp\left(-a^{\frac{1}{p}}\delta+\frac{\kappa(\eta)a^{\frac{2}{p}}}{4}
\right).
\end{align*}
The last inequality holds for arbitrary $a$ positive. Minimizing on
$a$ one finally obtains the desired estimate.
\end{proof}
\begin{prpstn}[Exponential tail estimates]\label{prpstn4}
If $Z$ is defined by $Z(t)=\int_0^tU(t-s)\xi(s)dW(s)$ such that
there exists $\eta$ positive such that
$\|\xi\|_{\xCn\left([0,T];\xHone(\R^d)\right)}^2\leq \eta$ a.s.,
then for any $p$ in $\tilde{\mathcal{A}}(d)$, $T$ and $\delta$
positive,
\begin{align*}
&\mathbb{P}\left(\|Z\|_{\xCn\left([0,T];\xHone(\R^d)\right)}\geq
\delta\right)\leq
3\exp\left(-\frac{\delta^2}{\kappa_1(\eta)}\right)\\
&\mathbb{P}\left(\|Z\|_{\xLn^{r(p)}\left(0,T;W^{1,p}(\R^d)\right)}\geq
\delta\right)\leq
c\exp\left(-\frac{\delta^2}{\kappa_2(\eta)}\right)\end{align*} where
$c=2e+\exp\left((2ek_0!)^{\frac{1}{k_0}}\right)$,
$k_0=2\vee\min\{k\in\N:2k\geq r(p)\}$
$$\kappa_1(\eta)=T4c\left(\infty\right)^2\|\Phi\|
_{\mathcal{L}_2^{0,s}}^2\eta,$$
$$\kappa_2(\eta)=\frac{8c\left(\frac{r(p)d}{2}\right)^2
T^{1-\frac{2}{r(p)}}(d+1)(d+p)\|\Phi\|_{\mathcal{L}_2^{0,s}}^2}
{1-\frac{4}{r(p)}}\eta.$$
\end{prpstn}
\begin{proof}
{\bf The first estimate.} We recall that in $\xHone(\R^d)$ we may
write, using the series expansion of the Wiener process and the fact
that $\left(U_t\right)_{t\in\R}$ is a unitary group, see \cite{EG},
$Z(t)=U(t)\int_0^tU(-s)\xi(s)dW(s)$. Since $U(-s)$ is an isometry,
one obtains that for every $s$ in $[0,T]$,
\begin{align*}
\|U(-s)\xi(s)\Phi\|_{\mathcal{L}_2(\L2,\xHone(\R^d))}&\leq
\|L_s\|_{\mathcal{L}_c\left(\xHn^s(\R^d,\R),\xHone(\R^d)\right)}
\|\Phi\|_{\mathcal{L}_2^{0,s}}\\
&\leq c\left(\infty\right)\|L_s\|_{\mathcal{L}_c\left(\xW^{1,\infty}
(\R^d,\R),\xHone(\R^d)\right)}
\|\Phi\|_{\mathcal{L}_2^{0,s}}\\
&\leq
c\left(\infty\right)\|\xi(s)\|_{\xHone(\R^d)}\|\Phi\|_{\mathcal{L}_2^{0,s}}
\end{align*}
where $L$ is such that $L_su=\xi(s)u$. Consequently, we obtain that
$$\int_0^T\|U(-s)\xi(s)\Phi\|_{\mathcal{L}_2^{0,1}}^2ds\leq c
\left(\infty\right)^2\|\Phi\|_{\mathcal{L}_2^{0,s}}^2T\|\xi(s)\|
_{\xCn([0,T],\xHone(\R^d))}^2ds.$$ We conclude using Theorem 2.1 of
\cite{P2}. The result is indeed stated for operator integrands from
a Hilbert space $H$ into $H$ and it still holds when the operator
takes its value in another Hilbert space.\\
{\bf The second estimate.} From Markov's inequality it is enough to
show that
$$\mathbb{E}\left[\exp\left(\frac{1}{\kappa_2(\eta)}\left\|\int_0^{\cdot}
U(\cdot-s)\xi(s)dW(s)\right\|_{\xLn^{r(p)}(0,T;\xW^{1,p}(\R^d))}^2\right)
\right]\leq c.$$ For $k\geq k_0$, Jensen's inequality along with
Fubini's theorem, Lemma \ref{l1}, the change of variables and the
integration by parts formulae give that
\begin{center}
\begin{tabular}{l}
$\mathbb{E}\left[\left(\frac{1}{\sqrt{\kappa_2(\eta)}}\left\|\int_0^{\cdot}
U(\cdot-s)\xi(s)dW(s)\right\|_{\xLn^{r(p)}(0,T;W^{1,p}(\R^d))}\right)^{2k}\right]$\\
$\leq\frac{1}{T}\int_0^T\mathbb{E}\left[\left(\frac{T}{\left(\kappa_2(\eta)
\right)^{\frac{r(p)}{2}}}
\left\|\int_0^{t}U(t-s)\xi(s)dW(s)\right\|_{\xW^{1,p}(\R^d)}^{r(p)}
\right)^{\frac{2k}{r(p)}}\right]dt$\\
$\leq\frac{1}{T}\int_0^T\int_0^{\infty}\mathbb{P}\left(\left\|\int_0^{t}U(t-s)
\xi(s)dW(s)\right\|_{\xW^{1,p}(\R^d)}\geq
\left(\frac{\left(\kappa_2(\eta)\right)^{\frac{r(p)}{2}}}{T}
\right)^{\frac{1}{r(p)}}u^{\frac{1}{2k}}\right) dudt$\\
$\leq\frac{1}{T}\int_0^T\int_0^{\infty}e\exp\left(-
\frac{\kappa_2(\eta)}{T^{\frac{2}{r(p)}}}\frac{u^{\frac{1}{k}}}
{\kappa(\eta)}\right)dudt$\\
$\leq e\int_0^{\infty}\exp\left(-2u^{\frac{1}{k}}\right)du=
e\int_0^{\infty}kv^{k-1}\exp(-2v)dv=
2e\int_0^{\infty}v^k\exp(-2v)dv.$
\end{tabular}\end{center}
Thus, using Fubini's theorem one obtains that
\begin{center}
\begin{tabular}{l}
$\mathbb{E}\left[\left(\frac{1}{\sqrt{\kappa_2(\eta)}}
\left\|\int_0^{\cdot}U(\cdot-s)\xi(s)dW(s)\right\|_{\xLn^{r(p)}(0,T;W^{1,p}(\R^d))}
\right)^{2k_0}\right]$\\
$\leq k_0!\sum_{k\geq
k_0}\frac{1}{k!}\mathbb{E}\left[\left(\frac{1}{\sqrt{\kappa_2(\eta)}}
\left\|\int_0^{\cdot}U(\cdot-s)\xi(s)dW(s)\right\|_{\xLn^{r(p)}(0,T;W^{1,p}(\R^d))}
\right)^{2k}\right]$\\
$\leq k_0!\sum_{k\geq
k_0}\frac{1}{k!}2e\int_0^{\infty}v^k\exp(-2v)dv$\\
$\leq k_0!\sum_{k\in\N}\frac{1}{k!}2e\int_0^{\infty}v^k\exp(-2v)dv
=2e k_0!,$
\end{tabular}\end{center}
hence using H\"older's inequality one obtains
\begin{center}
\begin{tabular}{l}
$\sum_{k=0}^{k_0-1}\frac{1}{k!}\mathbb{E}\left[\left(\frac{1}{\sqrt{\kappa_2(\eta)}}
\left\|\int_0^{\cdot}U(\cdot-s)\xi(s)dW(s)\right\|_{\xLn^{r(p)}(0,T;W^{1,p}(\R^d))}
\right)^{2k}\right]$\\
$=\sum_{k=0}^{k_0-1}\frac{1}{\kappa_2(\eta)^{k}k!}\mathbb{E}\left[
\left\{\left(\left\|\int_0^{\cdot}U(\cdot-s)\xi(s)dW(s)
\right\|_{\xLn^{r(p)}(0,T;W^{1,p}(\R^d))}\right)^{2k_0}\right\}^{\frac{k}{k_0}}\right]$\end{tabular}\\
\begin{tabular}{l}
$\leq\sum_{k=0}^{k_0-1}\frac{1}{\kappa_2(\eta)^{k}k!}\mathbb{E}\left[\left(
\left\|\int_0^{\cdot}U(\cdot-s)\xi(s)dW(s)
\right\|_{\xLn^{r(p)}(0,T;W^{1,p}(\R^d))}\right)^{2k_0}\right]^{\frac{k}{k_0}}$\\
$\leq\sum_{k=0}^{k_0-1}\frac{\left[\left(2ek_0!\right)^{\frac{1}{k_0}}\right]^k}
{k!}\leq\exp\left((2ek_0!)^{\frac{1}{k_0}}\right).$
\end{tabular}\end{center}\vspace{0.05cm}

\noindent The end of the proof is now straightforward.
\end{proof}

\section{Proof of the uniform large deviation principle}\label{s4}

\subsection{Almost continuity of the It\^{o} map}\label{s41} The proof of the
uniform large deviation principle now relies on
\begin{prpstn}\label{prpstn5}
For every positive $a$, $R$ and $\rho$, $u_0$ in $\xHone(\R^d)$, $f$
in $C_a$, $T<\mathcal{T}\left(S(u_0,f)\right)$, $p$ in
$\mathcal{A}(d)$, there exists positive $\epsilon_0$, $\gamma$ and
$r$ such that for every $\epsilon$ in $(0,\epsilon_0]$ and
$\tilde{u}_0$ in $B_{\xHone(\R^d)}(u_0,r)$,
$$\epsilon\log\mathbb{P}\left(\left\|u^{\epsilon,\tilde{u}_0}-S(u_0,f)
\right\|_{X^{(T,p)}}\geq
\rho;\left\|\sqrt{\epsilon}W-f\right\|_{\xCn\left([0,T];\xHn^s(\R^d,\R)\right)}
<\gamma\right)\leq-R.$$
\end{prpstn}
\begin{proof} Take $u_0$ in $\xHone(\R^d)$, $f$ in $C_a$, $a$, $R$
and $\rho$ positive, $T<\mathcal{T}\left(S(u_0,f)\right)$ and $p$ in $\mathcal{A}(d)$.\\
{\bf Step 1: Change of measure to center the Wiener process around
$f$.} The function $f$ in $C_a$ is such that there exists $h$ in
$\xLtwo\left(0,T;\xLtwo(\R^d)\right)$ such that
$f(\cdot)=\int_0^{\cdot}\Phi h(s)ds$ and
$\frac{1}{2}\left\|h\right\|_{\xLtwo\left(0,T;\xLtwo(\R^d)\right)}^2\leq
a$. We denote by $W^{\epsilon}$ the process defined by
\begin{align*}W^{\epsilon}(t)=W(t)-\frac{1}{\sqrt{\epsilon}}\int_0^t\frac{\p
f}{\p s}ds &=W(t)-\frac{1}{\sqrt{\epsilon}}\int_0^t\Phi
h(s)ds\\
&=\Phi\left(W_c(t)-\frac{1} {\sqrt{\epsilon}}\int_0^t
h(s)ds\right).\end{align*} Since
$\mathbb{E}\left[\exp\left(\frac{1}{2}\int_0^T\|h(s)\|_{\xLtwo(\R^d)}^2
\right)ds\right]<+\infty,$ the Novikov condition is satisfied and
the Girsanov theorem gives that $W^{\epsilon}$ is a $\mu-$Wiener
process on $\xCn\left([0,T];\xHn^s(\R^d,\R)\right)$ under the
probability $\mathbb{P}^{\epsilon}$ defined by
$$\left.\frac{\xdif \mathbb{P}^{\epsilon}}{\xdif \mathbb{P}}\right|
_{\mathcal{F}_t}=\exp\left(\frac{1}{\sqrt{\epsilon}}\int_0^t
\left(h,dW_c(s)\right)_{\xLtwo(\R^d)}-\frac{1}{2\epsilon}\int_0^T
\|h(s)\|_{\xLtwo(\R^d)}^2ds\right).$$ Set
$U_{\epsilon}(t)=\exp\left(-\frac{1}{\sqrt{\epsilon}}\int_0^t
\left(h,dW_c(s)\right)_{\xLtwo(\R^d)}\right)$ and $\lambda$ such
that $a-\lambda<-R$ and denote by $A$ the event
$$\left\{\left\|u^{\epsilon,\tilde{u}_0}-S(u_0,f)\right\|_{X^{(T,p)}}\geq
\rho;\left\|\sqrt{\epsilon}W-f\right\|_{\xCn\left([0,T];\xHn^s(\R^d,\R)\right)}
<\gamma\right\},$$ then
\begin{align*}
\mathbb{P}(A)&=\mathbb{E}_{\mathbb{P}_{\epsilon}}\left\{\frac{\xdif
\mathbb{P}}{\xdif
\mathbb{P}_{\epsilon}}\indic_{A\cap\left\{U_{\epsilon}(T)\leq\exp
\left(\frac{\lambda}{\epsilon}\right)\right\}}\right\}+\mathbb{P}
\left(U_{\epsilon}(T)>\exp\left(\frac{\lambda}{\epsilon}\right)\right)\\
&\leq\mathbb{E}_{\mathbb{P}_{\epsilon}}\left\{\indic_A\exp\left(
\frac{\lambda}{\epsilon}+\frac{1}{2\epsilon}\int_0^T\|h(s)\|
_{\xLtwo(\R^d)}^2ds\right)\right\}+\exp\left(-\frac{\lambda}
{\epsilon}\right)\mathbb{E}\left(U_{\epsilon}(T)\right)\\
&\leq\exp\left(\frac{\lambda+a}{\epsilon}\right)\mathbb{P}
_{\epsilon}(A)+\exp\left(\frac{a-\lambda}{\epsilon}\right).
\end{align*}
The last inequality follows from the fact that
$$\left(\exp\left(-\frac{1}{\sqrt{\epsilon}}\int_0^t
\left(h,dW_c(s)\right)_{\xLtwo(\R^d)}-\frac{1}{2\epsilon}
\int_0^t\|h(s)\|_{\xLtwo(\R^d)}^2ds\right)\right)_{t\in[0,T]}$$ is a
uniformly integrable martingale. Finally we see that it is
sufficient to prove that there exists positive $\epsilon_0$,
$\gamma$ and $r$ such that for every $\epsilon$ in $(0,\epsilon_0]$
and $\tilde{u}_0$ in $B_{\xHone(\R^d)}(u_0,r)$,
$$\epsilon\log\mathbb{P}_{\epsilon}(A)\leq-R-\lambda-a,$$
or equivalently that
$$\epsilon\log\mathbb{P}_{\epsilon}\left(\left\|v^{\epsilon,\tilde{u}_0}-S(u_0,f)
\right\|_{X^{(T,p)}}\geq
\rho;\left\|\sqrt{\epsilon}W_{\epsilon}\right\|_{\xCn
\left([0,T];\xHn^s(\R^d,\R)\right)}<\gamma\right)\leq-R-\lambda-a,$$
where $v^{\epsilon,\tilde{u}_0}$ satisfies
$v^{\epsilon,\tilde{u}_0}(0)=\tilde{u}_0$ and
$$idv^{\epsilon,\tilde{u}_0}=\left(\Delta
v^{\epsilon,\tilde{u}_0}+\lambda|v^{\epsilon,
\tilde{u}_0}|^{2\sigma}v^{\epsilon,\tilde{u}_0}+\frac{\p f}{\p
t}v^{\epsilon,\tilde{u}_0}-\frac{i\epsilon}{2}F_{\phi}v^{\epsilon,
\tilde{u}_0}\right)dt
+\sqrt{\epsilon}v^{\epsilon,\tilde{u}_0}dW_{\epsilon}.$$ {\bf Step
2: Reduction to estimates of the stochastic
convolution.}\\
Remark, this is standard fact, that the unboundedness of the drift
and coefficient of the Wiener process is not a limitation since the
result of Proposition \ref{prpstn5} is local. A localisation
argument will therefore be used to overcome the apparent difficulty.
We replace $T$ by
$$\tau_{\rho}=\inf\left\{t:\ \|v^{\epsilon,\tilde{u}_0}-S(u_0,f)\|_{X^{(t,p)}}
\geq\rho\right\}\wedge T.$$ Since $T<\mathcal{T}(S(u_0,f))$,
$v^{\epsilon,\tilde{u}_0}$ necessarily satisfies
$$ \|v^{\epsilon,\tilde{u}_0}\|_{X^{(\tau_{\rho},p)}}\leq \rho+\|S(u_0,f)\|
_{X^{(\tau_{\rho},p)}}=D.$$ Using the estimates used in the proofs
of Proposition \ref{prpstn2} herein and of Theorem $4.1$ of
\cite{dBD1}, with a truncation in front of the nonlinearity of the
form $\theta\left(\frac{ \|S(u_0,f)\|_{X^{(s,p)}}}{D}\right)$ and
$\theta\left(\frac{
\|v^{\epsilon,\tilde{u}_0}\|_{X^{(s,p)}}}{D}\right)$, one may obtain
for some $t>0$ and $\nu\in\left(0,1-\frac{\sigma(d-2)}{2}\right)$
$$\left\|v^{\epsilon,\tilde{u}_0}-S(u_0,f)\right\|_{X^{(t\wedge\tau_{\rho},p)}}\leq
C\|\tilde{u}_0-u_0\|_{\H1}+\sqrt{\epsilon}\left\|\int_0^{\cdot}U(\cdot-s)
v^{\epsilon,\tilde{u}_0}(s)dW_{\epsilon}(s)\right\|_{X^{(t\wedge\tau_{\rho},p)}}$$
$$+C\left[(t\wedge\tau_{\rho})^{\nu}(D^{2\sigma})(1+D)+
(t\wedge\tau_{\rho})^{\frac{1}{2}-\frac{1}{r(p)}}\sqrt{a}+\epsilon
(t\wedge\tau_{\rho})^{1-\frac{2}{r(p)}}\right]
\left\|v^{\epsilon,\tilde{u}_0}-S(u_0,f)\right\|_{X^{(t\wedge\tau_{\rho},p)}}$$
$$+C\epsilon
(t\wedge\tau_{\rho})^{1-\frac{2}{r(p)}}\|S(u_0,f)\|_{X^{(t\wedge\tau_{\rho},p)}}.$$
Set $\epsilon\leq1$, then for $t=t^*$ small enough we obtain
\begin{align*}
\left\|v^{\epsilon,\tilde{u}_0}-S(u_0,f)\right\|_{X^{(t^*\wedge\tau_{\rho},p)}}
\leq &2\left(C\|\tilde{u}_0-u_0\|_{\H1}+C\epsilon D\right.\\
&\left.+\sqrt{\epsilon}\left\|\int_0^{\cdot}U(\cdot-s)
v^{\epsilon,\tilde{u}_0}(s)dW_{\epsilon}(s)\right\|_{X^{(\tau_{\rho},p)}}\right).
\end{align*}
Set
$N=\left\lfloor\frac{\tau_{\rho}}{t^*\wedge\tau_{\rho}}\right\rfloor,$
and for $i$ in $\left\{0,...,N\right\}$, $T_i=iT^*$ and $T_{N+1}=T$.
The previous bound also holds for
$\left\|v^{\epsilon,\tilde{u}_0}-S(u_0,f)\right\|_{X^{(T_i,T_{i+1},p)}}$
for every $i$ in $\left\{0,...,N\right\}$, replacing $\|y-x\|_{\H1}$
by
$\left\|v^{\epsilon,\tilde{u}_0}(T_i)-S(u_0,f)(T_i)\right\|_{\H1}$.\\
As for $i$ in $\left\{1,...,N\right\}$,
$\left\|v^{\epsilon,\tilde{u}_0}(T_i)-S(u_0,f)(T_i)\right\|_{\H1}
\leq\left\|v^{\epsilon,\tilde{u}_0}-S(u_0,f)\right\|_{X^{(T_{i-1},T_{i},p)}}$,
we obtain using the triangle inequality that
\begin{align*}
\left\|v^{\epsilon,\tilde{u}_0}-S(u_0,f)\right\|_{X^{(\tau_{\rho},p)}}\leq&2(N+1)\left(\sqrt{\epsilon}\left\|\int_0^{\cdot}U(\cdot-s)
v^{\epsilon,\tilde{u}_0}(s)dW_{\epsilon}(s)\right\|_{X^{(\tau_{\rho},p)}}+C\epsilon
D\right)\\
&+2C\sum_{i=1}^{N-1}\left\|v^{\epsilon,\tilde{u}_0}-S(u_0,f)
\right\|_{X^{(T_{i-1},T_{i},p)}}+2C\|u_0-\tilde{u}_0\|_{\H1}\\
\leq&2(N+1)\left(\sum_{i=0}^{N-1}(2C)^i\right)\left(\sqrt{\epsilon}\left\|
\int_0^{\cdot}U(\cdot-s)v^{\epsilon,\tilde{u}_0}(s)dW_{\epsilon}(s)
\right\|_{X^{(\tau_{\rho},p)}}\right.\\
&\left.+C\epsilon D\right)+(2C)^N\|u_0-\tilde{u}_0\|_{\H1}.
\end{align*}
We may suppose that $2C>1$ and thus it is enough to show that there
exists positive $\epsilon_0$, $\gamma$ and $r$ such that
$(2C)^Nr<\rho$ and for every $\epsilon$ in $(0,\epsilon_0]$ and
$\tilde{u}_0$ in $B_{\xHone(\R^d)}(u_0,r)$,\\
\begin{tabular}{r}
$\epsilon\log\mathbb{P}_{\epsilon}\left(\sqrt{\epsilon}
\left\|\int_0^{\cdot}U(\cdot-s)v^{\epsilon,\tilde{u}_0}(s)dW_{\epsilon}(s)
\right\|_{X^{(\tau_{\rho},p)}}+C\epsilon D\geq
\frac{(2C-1)(\rho-(2C)^Nr)}{2(N+1)\left((2C)^N-1\right)};\right.$\\
$\left.\left\|\sqrt{\epsilon}W_{\epsilon}\right\|_{\xCn
\left([0,T];\xHn^s(\R^d,\R)\right)}<\gamma\right)\leq-R-\lambda-a.$\end{tabular}\vspace{0.3cm}

\noindent {\bf Step 3: The case of the stochastic convolution.} We
now need to see that for fixed $u_0$ in $\xHone(\R^d)$, $f$ in
$C_a$, $a$, $R$ and $\rho$ positive,
$T<\mathcal{T}\left(S(u_0,f)\right)$ and $p$ in $\mathcal{A}(d)$,
there exists $\epsilon_0$, $\gamma$ and $r$ positive such that for
all $\epsilon$ in $(0,\epsilon_0]$ and $\tilde{u}_0$ in
$B_{\xHone(\R^d)}(u_0,r)$,
$$\epsilon\log\mathbb{P}_{\epsilon}\left(\sqrt{\epsilon}
\left\|\int_0^{\cdot}U(\cdot-s)v^{\epsilon,\tilde{u}_0}(s)dW_{\epsilon}(s)
\right\|_{X^{(\tau_{\rho},p)}}\geq
\rho;\left\|\sqrt{\epsilon}W_{\epsilon}\right\|_{\xCn
\left([0,\tau_{\rho}];\xHn^s(\R^d,\R)\right)}<\gamma\right)\leq-R.$$
We introduce an approximation of $v^{\epsilon,\tilde{u}_0}$ similar
to that of the proof of Proposition \ref{prpstn2}. That is for $n$
in $\mathbb{N}$ and $i$ in $\left\{0,...,n\right\}$, we set
$t_i=\frac{i\tau_{\rho}}{n}$ and define $v^{\epsilon,\tilde{u}_0,n}$
when $t_i\leq t<t_{i+1}$ by
$$v^{\epsilon,\tilde{u}_0,n}(t)=U(t-t_i)\left(v^{\epsilon,\tilde{u}_0}(t_i)\right).$$
For any $\delta$ positive we may write
$$\mathbb{P}_{\epsilon}\left(\sqrt{\epsilon}\left\|\int_0^{\cdot}
U(\cdot-s)v^{\epsilon,\tilde{u}_0}(s)dW_{\epsilon}(s)\right\|_{X^{(\tau_{\rho},p)}}\geq
\rho;\left\|\sqrt{\epsilon}W_{\epsilon}\right\|
_{\xCn([0,\tau_{\rho}];\xHn^s(\R^d,\R)}<\gamma\right)$$
\begin{tabular}{r r}
$\leq$ &
$\mathbb{P}_{\epsilon}\left(\sqrt{\epsilon}\left\|\int_0^{\cdot}U(\cdot-s)
\left(v^{\epsilon,\tilde{u}_0}(s)-v^{\epsilon,\tilde{u}_0,n}(s)\right)
dW_{\epsilon}(s)\right\|
_{X^{(\tau_{\rho},p)}}\geq\frac{\rho}{2};\right.$\\
&$\left.\left\|v^{\epsilon,\tilde{u}_0}-v^{\epsilon,\tilde{u}_0,n}
\right\|_{\xCn([0,\tau_{\rho}];\H1)}<\delta\right)$\end{tabular}\\
\begin{tabular}{r r}
\hspace{0.35cm}$+\mathbb{P}_{\epsilon}\left(\left\|v^{\epsilon,\tilde{u}_0}
-v^{\epsilon,\tilde{u}_0,n}
\right\|_{\xCn([0,\tau_{\rho}];\H1)}\geq\delta\right)$\end{tabular}\\
\begin{tabular}{r r}
&$+\mathbb{P}_{\epsilon}\left(\sqrt{\epsilon}\left\|\int_0^{\cdot}
U(\cdot-s)v^{\epsilon,\tilde{u}_0,n}(s)dW_{\epsilon}(s)\right\|_{X^{(\tau_{\rho},p)}}
\geq \frac{\rho}{2};\left\|\sqrt{\epsilon}W_{\epsilon}\right\|_{\xCn
\left([0,\tau_{\rho}];\xHn^s(\R^d,\R)\right)}<\gamma;\right.$\\
&$\left.\left\|v^{\epsilon,\tilde{u}_0}-v^{\epsilon,\tilde{u}_0,n}\right\|
_{\xCn([0,\tau_{\rho}];\H1)}<\delta\right)$
\end{tabular}\vspace{0.3cm}

\noindent {\bf Bound for the first term.} From Proposition
\ref{prpstn4},
$C\exp\left(-\frac{\rho^2}{4\epsilon\left(\kappa_1(\delta^2)
\vee\kappa_2(\delta^2)\right)}\right)$ is an upper bound for the
first  term. Thus for any $\epsilon$ positive and $\delta$ small
enough\\
\begin{tabular}{r}
$\epsilon\log\mathbb{P}_{\epsilon}\left(\sqrt{\epsilon}
\left\|\int_0^{\cdot}U(\cdot-s)\left(v^{\epsilon,\tilde{u}_0}(s)-
v^{\epsilon,\tilde{u}_0,n}(s)\right)dW_{\epsilon}(s)\right\|
_{X^{(\tau_{\rho},p)}}\geq \frac{\rho}{2};\right.$\\
$\left.\left\|v^{\epsilon,\tilde{u}_0}-v^{\epsilon,\tilde{u}_0,n}
\right\|_{\xCn([0,\tau_{\rho}];\H1)}<\delta\right)$\end{tabular}\vspace{0.3cm}

\noindent is less than $-R-1$.\\
{\bf Bound for the second term.} The second term
$\mathbb{P}_{\epsilon}\left(\left\|v^{\epsilon,\tilde{u}_0}-
v^{\epsilon,\tilde{u}_0,n}\right\|_{\xCn([0,\tau_{\rho}];\H1)}
\geq\delta\right)$ is bounded by the sum of
$$\mathbb{P}_{\epsilon}\left(\sup_{i\in\left\{0,...,n-1\right\}}
\sqrt{\epsilon}\left\|\int_{t_i}^{\cdot}U(t-s)v^{\epsilon,\tilde{u}_0}(s)
dW_{\epsilon}(s)\right\|_{\xCn([t_i,t_{i+1}];\H1)}\geq\frac{\delta}{2}\right)$$
which from Proposition \ref{prpstn4} is less than
$3n\exp\left(-\frac{C\delta^2n}{\tau_{\rho}D^2}\right)$ and
$$\mathbb{P}_{\epsilon}\left(\sup_{i\in\left\{0,...,n-1\right\}}
\left\|\int_{t_i}^{\cdot}U(t-s)\left[\lambda\left|
v^{\epsilon,\tilde{u}_0}(s) \right|^{2\sigma}+\frac{\p f}{\p
s}(s)-\frac{i\epsilon}{2}F_{\Phi}\right]
v^{\epsilon,\tilde{u}_0}(s)ds\right\|_{\xCn([t_i,t_{i+1}];\H1)}\geq
\frac{\delta}{2}\right).$$ The second term is equal to zero for $n$
large enough. Indeed, with calculations similar to that of the proof
of Theorem $4.1$ of \cite{dBD1} and of Proposition \ref{prpstn2}, we
obtain that for $\epsilon<1$,
$$\sup_{i\in\left\{0,...,n-1\right\}}\left\|\int_{t_i}^{\cdot}
U(t-s)\left[\lambda\left|
v^{\epsilon,\tilde{u}_0}(s)\right|^{2\sigma}+\frac{\p f}{\p s}(s)
-\frac{i\epsilon}{2}F_{\Phi}\right]
v^{\epsilon,\tilde{u}_0}(s)ds\right\|_{\xCn([t_i,t_{i+1}];\H1)}$$
$$\leq
C\left[\left(\frac{\tau_{\rho}}{n}\right)^{\nu}D^{2\sigma+1}+
\left(\frac{\tau_{\rho}}{n}\right)^{\frac{1}{2}-\frac{1}{r(p)}}D\sqrt{a}
+\frac{1}{2}\left(\frac{\tau_{\rho}}{n}\right)^{1-\frac{2}{r(p)}}D\right].$$
Finally, for every $\delta$ positive and $0<\epsilon<1$, for $n$
large enough
$$\epsilon\log\mathbb{P}_{\epsilon}\left(\left\|v^{\epsilon,\tilde{u}_0}
-v^{\epsilon,\tilde{u}_0,n}
\right\|_{\xCn([0,\tau_{\rho}];\H1)}\geq\delta\right)<-R-1.$$ Thus
choosing first $\delta$ small enough and then $n$ large enough we
obtain that for any $0<\epsilon<\frac{1}{2\log(2)}$,\\
\begin{tabular}{r}
$\epsilon\log\left\{\mathbb{P}_{\epsilon}\left(\sqrt{\epsilon}
\left\|\int_0^{\cdot}U(\cdot-s)
\left(v^{\epsilon,\tilde{u}_0}(s)-v^{\epsilon,\tilde{u}_0,n}(s)\right)
dW_{\epsilon}(s)\right\|_{X^{(\tau_{\rho},p)}}\geq\frac{\rho}{2};\right.\right.$\\
$\left.\left.\left\|v^{\epsilon,\tilde{u}_0}-v^{\epsilon,\tilde{u}_0,n}
\right\|_{\xCn([0,\tau_{\rho}];\H1)}<\delta\right)\right.$\end{tabular}

\begin{tabular}{r}
 \hspace{0.8cm} $\left.+\mathbb{P}_{\epsilon}\left(\left\|v^{\epsilon,\tilde{u}_0}
-v^{\epsilon,\tilde{u}_0,n}
\right\|_{\xCn([0,\tau_{\rho}];\H1)}\geq\delta\right)\right\}
\leq-R-\frac{1}{2}$
\end{tabular}\vspace{0.3cm}

\noindent {\bf Bound for the third term.} Fix $\delta$ as above.
Remark first that the last condition
$$\left\|v^{\epsilon,\tilde{u}_0}-v^{\epsilon,\tilde{u}_0,n}\right\|
_{\xCn([0,\tau_{\rho}];\H1)}<\delta$$ implies that
$\left\|v^{\epsilon,\tilde{u}_0,n}\right\|_{\xCn([0,\tau_{\rho}];\H1)}<D+\delta$.
Denote by $$\underline{t}=\max\left\{t_i:t_i\leq t,\
i\in\{0,...,n\}\right\}$$and by $E$ the event
$$E=\left\{\left\|\sqrt{\epsilon}W_{\epsilon}\right\|
_{\xCn([0,\tau_{\rho}];\xHn^s(\R^d,\R))}<\gamma\ ;\
\left\|v^{\epsilon,\tilde{u}_0}-v^{\epsilon,\tilde{u}_0,n}
\right\|_{\xCn([0,\tau_{\rho}];\H1)}<\delta\right\}.$$ As
$p<\frac{2(3d-1)}{3(d-1)}$ there exists $p<\tilde{p}<\frac{2d}{d-1}$
and $\eta$ positive such that
$1-\frac{p-2}{\tilde{p}-2}\left(1+\frac{2}{r(\tilde{p})}+\eta\right)$
is positive. Thus from H\"older's inequality, setting
$\theta=\frac{\tilde{p}-p}{\tilde{p}-2}$, the third term is bounded
above by\\
\begin{tabular}{r r}
   & $\mathbb{P}_{\epsilon}\left( \sqrt{\epsilon}\left\|
\int_{\underline{\cdot}}^{\cdot}U(\cdot-s)
v^{\epsilon,\tilde{u}_0,n}(s)dW_{\epsilon}(s)\right\|_{\xLn^{r(\tilde{p})}
\left(0,\tau_{\rho};\xW^{1,\tilde{p}}(\R^d)\right)}\geq
n^{\left(1+\frac{2}{r(\tilde{p})}+\eta\right)\frac{1}{2}};\right.$\\
  & $\left.\left\|v^{\epsilon,\tilde{u}_0}-v^{\epsilon,\tilde{u}_0,n}
\right\|_{\xCn([0,\tau_{\rho}];\H1)}<\delta\right)$\\
  $+$ & $\mathbb{P}_{\epsilon}\left(\sqrt{\epsilon}\left\|
\int_{\underline{\cdot}}^{\cdot}U(\cdot-s)
v^{\epsilon,\tilde{u}_0,n}(s)dW_{\epsilon}(s)
\right\|_{\xCn\left([0,\tau_{\rho}];\H1\right)}^{\theta}
n^{\left(1+\frac{2}{r(\tilde{\rho})}+\eta\right)
\frac{1-\theta}{2}}\geq\frac{\rho}{8};\right.$\\
   & $\left.\left\|v^{\epsilon,\tilde{u}_0}-
v^{\epsilon,\tilde{u}_0,n}\right\|_{\xCn([0,\tau_{\rho}];\H1)}<\delta\right)$
\end{tabular}\\
\begin{tabular}{r r}
  $+$ & $\mathbb{P}_{\epsilon}\left(\sqrt{\epsilon}\left\|
\int_{\underline{\cdot}}^{\cdot}U(\cdot-s)
v^{\epsilon,\tilde{u}_0,n}(s)dW_{\epsilon}(s)
\right\|_{\xCn\left([0,\tau_{\rho}];\H1\right)}\geq\frac{\rho}{8};\right.$\\
   & $\left.\left\|v^{\epsilon,\tilde{u}_0}-
v^{\epsilon,\tilde{u}_0,n}\right\|_{\xCn([0,\tau_{\rho}];\H1)}<\delta\right)$
\end{tabular}\\
\begin{tabular}{r r}
  $+$ & $\mathbb{P}_{\epsilon}\left(\sqrt{\epsilon}\left\|\int_0^{\underline{\cdot}}
U(\cdot-s)v^{\epsilon,\tilde{u}_0,n}(s)dW_{\epsilon}(s)\right\|_{\xLn^{r(p)}
(0,\tau_{\rho};W^{1,p}(\R^d))}\geq \frac{\rho}{8};E\right)$
\end{tabular}\\
\begin{tabular}{r r}
  $+$ & $\mathbb{P}_{\epsilon}\left(\sqrt{\epsilon}\left\|
\int_0^{\underline{\cdot}}U(\cdot-s)v^{\epsilon,\tilde{u}_0,n}(s)
dW_{\epsilon}(s)\right\|_{\xCn([0,\tau_{\rho}];\H1)}\geq
\frac{\rho}{8};E\right).$\\
 $+$ & $\mathbb{P}_{\epsilon}\left(\sqrt{\epsilon}\left\|
\int_0^{\underline{\cdot}}U(\cdot-s)v^{\epsilon,\tilde{u}_0,n}(s)
dW_{\epsilon}(s)\right\|_{\xCn([0,\tau_{\rho}];\H1)}\geq
\frac{\rho}{8};E\right).$
\end{tabular}\vspace{0.3cm}

\noindent The first probability is bounded from above by\\
\begin{tabular}{r}
$\sum_{i=0}^{n-1}\mathbb{P}_{\epsilon}\left(\sqrt{\epsilon}
\left\|\int_{t_i}^{t_{i+1}}U(\cdot-s)
v^{\epsilon,\tilde{u}_0,n}(s)dW_{\epsilon}(s)
\right\|_{\xLn^{r(\tilde{p})}(t_i,t_{i+1};W^{1,\tilde{p}}(\R^d))}\geq
n^{\left(-1+\frac{2}{r(\tilde{p})}+\eta\right)\frac{1}{2}};\right.$\\
$\left.\left\|v^{\epsilon,\tilde{u}_0}-v^{\epsilon,\tilde{u}_0,n}
\right\|_{\xCn([0,\tau_{\rho}];\H1)}<\delta\right)$
\end{tabular}\vspace{0.3cm}

\noindent which from Proposition \ref{prpstn4} is less than
$$nC\exp\left(-\frac{n^{-1+\frac{2}{r(\tilde{p})}
+\eta}}{C\left(\frac{\tau_{\rho}}{n}\right)^{1-
\frac{2}{r(\tilde{p})}}\epsilon(\delta+D)^2}\right).$$ Thus for $n$
large enough, for any $\epsilon$ and $\delta$ positive,\\
\begin{tabular}{r}
$\epsilon\log\mathbb{P}_{\epsilon}\left(\sqrt{\epsilon}\left\|
\int_{\underline{\cdot}}^{\cdot}U(\cdot-s)
v^{\epsilon,\tilde{u}_0,n}(s)dW_{\epsilon}(s)\right\|_{\xLn^{r(\tilde{p})}
\left(0,\tau_{\rho};\xW^{1,\tilde{p}}(\R^d)\right)}
\geq n^{\left(1+\frac{2}{r(\tilde{p})}+\eta\right)\frac{1}{2}};\right.$\\
$\left.\left\|v^{\epsilon,\tilde{u}_0}-v^{\epsilon,\tilde{u}_0,n}
\right\|_{\xCn([0,\tau_{\rho}];\H1)}<\delta\right)<-R-1.$
\end{tabular}\vspace{0.3cm}

\noindent The first exponential tail estimate of Proposition
\ref{prpstn4} gives that the second probability is less than
$$nC\exp\left(-\frac{\rho^2}{Cn^{\left(1+\frac{2}{r(\tilde{p})}+\eta\right)(1-\theta)}
\left(\frac{\tau_{\rho}}{n}\right)\epsilon(\delta+D)^2}\right),$$
and, from the choice of $\tilde{p}$, for $n$ large enough, for any
$\epsilon$ and $\delta$ positive,\\
\begin{tabular}{r}
$\epsilon\log\mathbb{P}_{\epsilon}\left(\sqrt{\epsilon}\left\|
\int_{\underline{\cdot}}^{\cdot}U(\cdot-s)
v^{\epsilon,\tilde{u}_0,n}(s)dW_{\epsilon}(s)
\right\|_{\xCn\left([0,\tau_{\rho}];\H1\right)}^{\theta}
n^{\left(1+\frac{2}{r(\tilde{\rho})}+\eta\right)
\frac{1-\theta}{2}}\geq\frac{\rho}{3};\right.$\\
$\left.\left\|v^{\epsilon,\tilde{u}_0}-v^{\epsilon,
\tilde{u}_0,n}\right\|_{\xCn([0,\tau_{\rho}];\H1)}<\delta\right)<-R-1.$
\end{tabular}\vspace{0.3cm}

\noindent The same holds for the third probability more clearly.\\
The decay estimates of Section \ref{s22} along with H\"older's
inequality give that the mapping
$$w\mapsto U(t-t_j)v^{\epsilon,\tilde{u}_0}(t_j)w$$
from $\xHn^s(\R^d,\R)$ to $\xW^{1,p}(\R^d)$ is continuous. Thus, we
may write
\begin{center}
\begin{tabular}{l}
$\left\|\int_0^{\underline{\cdot}}U(\cdot-s)v^{\epsilon,
\tilde{u}_0,n}(s)dW_{\epsilon}(s)
\right\|_{\xLn^{r(p)}(0,T;W^{1,p}(\R^d))}$\\
$=\left\|\sum_{i=1}^{n-1}\indic_{t_i\leq
t<t_{i+1}}\sum_{j=0}^{i-1}\int_{t_j}^{t_{j+1}}U(t-t_j)v^{\epsilon,\tilde{u}_0}(t_j)
dW_{\epsilon}(s)\right\|_{\xLn^{r(p)}(0,T;W^{1,p}(\R^d))}$\\
$\leq\sum_{i=1}^{n-1}\sum_{j=0}^{i-1}\left\|\int_{t_j}^{t_{j+1}}U(t-t_j)
v^{\epsilon,\tilde{u}_0}(t_j)dW_{\epsilon}(s)\right\|_{\xLn^{r(p)}(t_i,t_{i+1};
W^{1,p}(\R^d))}$\\
$\leq C\left(\frac{(n-1)(n-2)}{2}\right)\left(\frac{\tau_{\rho}}{n}
\right)^{-\frac{2}{r(p)}}D\gamma,$
\end{tabular}\end{center}\vspace{0.05cm}

and obtain that, for any $n$ in $\N$, for
$\gamma$ small enough the fourth probability is equal to zero.\\
Similarly we write, using the continuity of the group and H\"older's
inequality,\\
\begin{center}
\begin{tabular}{l}
$\left\|\int_0^{\underline{\cdot}}U(\cdot-s)v^{\epsilon,\tilde{u}_0,n}
(s)dW_{\epsilon}(s)
\right\|_{\xCn(0,T;\H1)}$\\
$=\max_{i=1,...,n-1}\left\|\sum_{j=0}^{i-1}\int_{t_j}^{t_{j+1}}U(t_i-t_j)
v^{\epsilon,\tilde{u}_0}(t_j)dW_{\epsilon}(s)\right\|_{\H1}$\\
$\leq\sum_{j=0}^{n-1}\|v^{\epsilon,\tilde{u}_0}(t_j)\|_{\H1}\|W_{\epsilon}(t_{j+1})
-W_{\epsilon}(t_j)\|_{\xHn^s(\R^d,\R)}$\\
$\leq 2nD\gamma.$
\end{tabular}\end{center}\vspace{0.05cm}

Thus, for any $n$ in $\N$ , for
$\gamma$ small enough the fifth probability is equal to zero.\\
Finally, when $\delta$ is fixed, for $n$ large enough and a
particular choice of $\gamma$ depending on $n$ and $\delta$, we
obtain that for any $0<\epsilon<\frac{1}{2\log(2)}$,\\
\begin{tabular}{r}
$\epsilon\log\mathbb{P}_{\epsilon}
\left(\sqrt{\epsilon}\left\|\int_0^{\cdot}
U(\cdot-s)v^{\epsilon,\tilde{u}_0,n}(s)dW_{\epsilon}(s)
\right\|_{X^{(\tau_{\rho},p)}}\geq
\frac{\rho}{2};\left\|\sqrt{\epsilon}W_{\epsilon}\right\|_{\xCn
\left([0,\tau_{\rho}];\xHn^s(\R^d,\R)\right)}<\gamma;\right.$\\
$\left.\left\|v^{\epsilon,\tilde{u}_0}-v^{\epsilon,\tilde{u}_0,n}\right\|
_{\xCn([0,\tau_{\rho}];\H1)}<\delta\right)\leq
-R-\frac{1}{2}.$\end{tabular}\vspace{0.3cm}

\noindent We have now proved Step 3 and thus Proposition
\ref{prpstn5}.
\end{proof}
\begin{rmrk}
In Step $2$ we did not use the Gronwall inequality, which is often
used in that case. Instead we split the norm in many parts because
in the case of the Schr\"odinger group it is needed, in order to use
the integrability property governed by $ii/$ of the Strichartz
estimates, to keep the convolution with the group.
\end{rmrk}
\begin{rmrk}
Revisiting Step 2 we may see that a uniform LDP holds with the same
good rate function and same speed for the laws of the paths of the
solutions for an equation with an additional term of the form
$f\left(u^{\epsilon,u_0},\epsilon,t,x\right)$ in the drift. It is
needed that there exists $(s,\rho)$ conjugate exponents of an
admissible pair $(r(q),q)$ such that for every positive $T$ such
that $\left\|\psi\right\|_{X^{(T,p)}}<+\infty$,
$\left\|f\left(\psi,\epsilon,.,*\right)\right\|_{\xLn^{s}\left(0,T;\xW^{1,\rho}(\R^d)\right)}$
is bounded and goes to zero as $\epsilon$ goes to zero. This extra
term may be an external potential accounting for damping or
amplification going to zero along with the noise intensity. In
\cite{CGJRY,GKR} for example, extra terms are added to the equation
to account for nonlinear, respectively linear, damping. In
\cite{FKLT} a small amplification added to compensate for loss due
to a small amplitude modulation is considered. The amplitude
modulation makes sense considering initial data in spaces of
spatially localized functions used when studying the blow-up of
deterministic and stochastic NLS equations.\end{rmrk}
\begin{rmrk}
The proof of local existence and uniqueness of the solutions of the
stochastic NLS equation in $\xHone(\R^d)$ holds with more general
nonlinearities and when the noise enters nonlinearly under some
Lipschitz assumptions. We may adapt our proof of the uniform LDP to
cover those NLS equations.
\end{rmrk}
\subsection{End of the proof}\label{s42} We prove hereafter how the almost
continuity along with Proposition \ref{prpstn1} and \ref{prpstn2}
allow to prove the uniform LDP for the laws of the solutions in the
space $\mathcal{E}_{\infty}$ when the noise goes to zero.\\
Suppose that $I^{u_0}$ is the rate function of the LDP then from
Proposition \ref{prpstn2}, since its level sets are the direct image
by $S(u_0,\cdot)$ of the level sets $C_a$ which are compact, it is a
good rate function.\\
The set $A$ is a Borel set of $\mathcal{E}_{\infty}$ and $u_0$ is
some initial datum in $\H1$.\\
\noindent{\bf An upper bound.} In the case where
$\inf_{w\in\overline{A}}I^{u_0}(w)=0$ there is nothing to prove.
Otherwise, take $0<a<\inf_{w\in\overline{A}}I^{u_0}(w)$ and $R>a$.
Suppose that $f$ is such that $I^W(f)\leq a$, then
$$I^{u_0}(S(u_0,f))\leq a<\inf_{w\in\overline{A}}I^{u_0}(w),$$
thus $S(u_0,f)\notin \overline{A}$ and there exists an elementary
neighborhood of $S(u_0,f)$ of the form
$$V_{u_0,f}=\left\{v\in\mathcal{E}_{\infty}:\ \mathcal{T}(v)>T\ \mbox{and}\
\|v-S(u_0,f)\|_{X^{(T,p)}}<\rho_{u_0,f}\right\}$$ for some $(p,T)$
such that $V_{u_0,f}\subset\overline{A}^{c}$. Also, from Proposition
\ref{prpstn4}, there exists $\epsilon_{u_0,f}$, $\gamma_{u_0,f}$ and
$r_{u_0,f}$ positive such that for every
$\epsilon\leq\epsilon_{u_0,f}$ and $\tilde{u}_0$ in
$B_{\H1}(u_0,r_{u_0,f})$,
$$\epsilon\log\mathbb{P}\left(\left\|u^{\epsilon,\tilde{u}_0}-S(u_0,f)
\right\|_{X^{(T,p)}}\geq
\rho_{u_0,f};\left\|\sqrt{\epsilon}W-f\right\|_{\xCn
\left([0,T];\xHn^s(\R^d,\R)\right)}<\gamma_{u_0,f}\right)\leq-R.$$
Let denote by $O_{u_0,f}$ the set
$$O_{u_0,f}=B_{\xCn\left([0,T];\xHn^s(\R^d,\R)
\right)}(f,\gamma_{u_0,f}).$$ The family
$\left(O_{u_0,f}\right)_{f\in C_a}$ is a covering by open sets of
the compact set $C_a$, thus there exists a finite sub-covering of
the form $\bigcup_{i=1}^NO_{u_0,f_i}.$ We can now write
\begin{align*}
\mathbb{P}\left(u^{\epsilon,\tilde{u}_0}\in A\right)\leq&
\mathbb{P}\left(\left\{u^{\epsilon,\tilde{u}_0}\in
A\right\}\cap\left\{\sqrt{\epsilon}W\in\bigcup_{i=1}^N
O_{u_0,f_i}\right\}\right)+\mathbb{P}\left(\sqrt{\epsilon}W
\notin\bigcup_{i=1}^NO_{u_0,f_i}\right)\\
\leq &\sum_{i=1}^N\mathbb{P}\left(\left\{u^{\epsilon,\tilde{u}_0}\in
A \right\}\cap\left\{\sqrt{\epsilon}W\in
O_{u_0,f_i}\right\}\right)+\mathbb{P}\left(\sqrt{\epsilon}W\notin
C_a\right)\\
\leq
&\sum_{i=1}^N\mathbb{P}\left(\left\{u^{\epsilon,\tilde{u}_0}\notin
V_{u_0,f_i} \right\}\cap\left\{\sqrt{\epsilon} W\in
O_{u_0,f_i}\right\}\right)+\exp\left(-\frac{a}{\epsilon}\right),
\end{align*}
for $\epsilon\leq\epsilon_0$ for some $\epsilon_0$ positive. Thus
for $\epsilon\leq\epsilon_0\wedge\left(\bigwedge_{i=1}^m
\epsilon_{u_0,f_i}\right)$ we obtain for $\tilde{u}_0$ in
$B_{\H1}(u_0,r_{u_0})$ where $r_{u_0}=\bigwedge_{i=1}^mr_{u_0,f_i}$,
$$\mathbb{P}\left(u^{\epsilon,\tilde{u}_0}\in A\right)\leq\
N\exp\left(-\frac{R}{\epsilon}\right)+\exp\left(-\frac{a}{\epsilon}\right),$$
and $$\epsilon\log\mathbb{P}\left(u^{\epsilon,\tilde{u}_0}\in
A\right)\leq\epsilon\log2+\left(\epsilon\log N-R\right)\vee(-a).$$
Finally, there exists $r_{u_0}$ such that for any $\tilde{u}_0$ in
$B_{\H1}(u_0,r_{u_0})$,
$$\overline{\lim}_{\epsilon\rightarrow0}\epsilon\log\mathbb{P}
\left(u^{\epsilon,\tilde{u}_0}\in A\right)\leq-a.$$ Since $a$ is
arbitrary, we obtain,
$$\overline{\lim}_{\epsilon\rightarrow0,\tilde{u}_0\rightarrow u_0}
\epsilon\log\mathbb{P}\left(u^{\epsilon,\tilde{u}_0}\in
A\right)\leq-\inf_{w\in\overline{A}}I^{u_0}(w).$$ {\bf A lower
bound.} Suppose that $\inf_{w\in Int(A)}I^{u_0}(w)<+\infty$,
otherwise there is nothing to prove, and take
$w$ in $Int(A)$ such that $I^{u_0}(w)<+\infty$.\\
The continuity of $S(u_0,\cdot)$ along with the compactness the
level set $C_{I^{u_0}(w)+1}$ give that there exists $f$ such that
$w=S(u_0,f)$ and $I^{u_0}(w)=I^W(f)$. Take $V_{u_0,f}$ an elementary
neighborhood of $S(u_0,f)$ included in $A$ and $O_{u_0,f}$ defined
as previously, $\eta$ positive and $R>I^{u_0}(w)+\eta$. We obtain
\begin{align*}
\exp\left(-\frac{R-\eta}{\epsilon}\right)&\leq\exp
\left(-\frac{I^W(f)}{\epsilon}\right)\\
&\leq\mathbb{P}\left(\sqrt{\epsilon}W\in
O_{u_0,f}\right)\\
&\leq\mathbb{P}\left(\left\{u^{\epsilon,\tilde{u}_0}\notin
V_{u_0,f}\right\}\cap\left\{\sqrt{\epsilon}W\in
O_{u_0,f}\right\}\right)+\mathbb{P}\left(u^{\epsilon,\tilde{u}_0}\in
A\right).
\end{align*}
Thus there exists $r_{u_0}$ and $\epsilon_0$ positive such that for
all $\tilde{u}_0$ in $B_{\H1}(u_0,r_{u_0})$ and
$\epsilon\leq\epsilon_0$,
$$-R+\eta\leq\epsilon\log 2+\left(\epsilon\log\mathbb{P}
\left(u^{\epsilon,\tilde{u}_0}\in A\right)\right)\vee(-R)$$ and
there exists $\epsilon_1\leq\epsilon_0$ such that for all
$\epsilon\leq\epsilon_1,$
$$-I^{u_0}(w)\leq\epsilon\log2+\epsilon\log\mathbb{P}\left(u^{\epsilon,\tilde{u}_0}\in
A\right).$$ As a consequence, we obtain that for every $u_0$ in
$\H1$, there exists $r_{u_0}$ positive such that for every
$\tilde{u}_0$ in $B_{\H1}(u_0,r_{u_0})$,
$$\underline{\lim}_{\epsilon\rightarrow0}\epsilon
\log\mathbb{P}\left(u^{\epsilon,\tilde{u}_0}\in
A\right)\geq-I^{u_0}(w)$$ and
$$\underline{\lim}_{\epsilon\rightarrow0}\epsilon
\log\mathbb{P}\left(u^{\epsilon,\tilde{u}_0}\in
A\right)\geq-\inf_{w\in Int(A)}I^{u_0}(w)$$ since $w$ in $Int(A)$ is
arbitrary.\\
{\bf Uniformity with respect to initial data in compact sets.} The
uniform LDP follows from the previous bounds along with Corollary
$5.6.15$ of \cite{DZ}. We now give the proof of the lower bound
which is not written in the previous reference. Let $K$ be a compact
set in $\H1$. Suppose that $\sup_{\tilde{u}_0\in K}\inf_{w\in
Int(A)}I^{\tilde{u}_0}(w)<+\infty$, otherwise there is nothing to
prove and take $\delta$ positive. For any $u_0$ in $K$, there exists
$r_{u_0}$ positive such that for every $\epsilon\leq\epsilon_{u_0}$,
$$\epsilon\log\inf_{\tilde{u}_0\in
B_{\H1}(u_0,r_{u_0})}\mathbb{P}\left(u^{\epsilon,\tilde{u}_0}\in
A\right)\geq -\sup_{\tilde{u}_0\in K}\inf_{w\in
Int(A)}I^{\tilde{u}_0}(w)-\delta.$$ The set of balls
$B_{\H1}(u_0,r_{u_0})$ is a covering of $K$ by open sets thus there
exists a sub-covering of $K$ of the form
$\bigcup_{i=1}^mB_{\H1}(u_0^i,r_{u_0^i})$ and for $\epsilon\leq
\bigwedge_{i=1}^m\epsilon_{u_0^i}$,
\begin{align*}
\epsilon\log\inf_{\tilde{u}_0\in
K}\mathbb{P}\left(u^{\epsilon,\tilde{u}_0}\in
A\right)&\geq\epsilon\log\inf_{\tilde{u}_0\in
\bigcup_{i=1}^mB_{\H1}(u_0^i,r_{u_0^i})}\mathbb{P}\left(u^{\epsilon,\tilde{u}_0}\in
A\right)\\
&\geq -\sup_{\tilde{u}_0\in K}\inf_{w\in
Int(A)}I^{\tilde{u}_0}(w)-\delta,\end{align*} conclusion follows
since $\delta$ is arbitrary.
\section{Applications to the blow-up time}\label{s5}
In this section the equation with a focusing nonlinearity is
considered. In this case, it is known that some solutions of the
deterministic equation blow up in finite time for critical or
supercritical nonlinearities. If $B$ is a Borel set of
$[0,+\infty]$,
$$\mathbb{P}\left(\mathcal{T}\left(u^{\epsilon,u_0}\right)\in
B\right)=\mu^{u^{\epsilon,u_0}}\left(\mathcal{T}^{-1}(B)\right).$$
Thus the uniform LDP for the family
$\left(\mu^{u^{\epsilon,u_0}}\right)_{\epsilon>0}$ gives that for
$K\subset\subset\H1$,
$$-\sup_{u_0\in K}\inf_{u\in
Int\left(\mathcal{T}^{-1}(B)\right)}I^{u_0}(u)\leq
\underline{\lim}_{\epsilon\rightarrow0}\epsilon \log\inf_{u_0\in
K}\mathbb{P}\left(\mathcal{T}(u^{\epsilon,u_0})\in B\right)$$and
that
$$\overline{\lim}_{\epsilon\rightarrow0}\epsilon
\log\sup_{u_0\in K}\mathbb{P}\left(\mathcal{T}(u^{\epsilon,u_0})\in
B\right)\leq-\inf_{u\in \overline{\mathcal{T}^{-1}(B)},u_0\in
K}I^{u_0}(u).$$ Since $\mathcal{T}$ is lower semicontinuous the sets
$(T,+\infty]$ and $[0,T]$ are particularly interesting. We recall,
see \cite{EG} for more details, that for every $T>0$,
$\overline{\mathcal{T}^{-1}((T,+\infty])}=\mathcal{E}_{\infty}$ and
$Int\left(\mathcal{T}^{-1}([0,T])\right)=\emptyset$. Thus for the
two types of sets, at least one bound is trivial. Considering the
approximate blow-up time allows us to obtain two interesting bounds
and to treat intervals of the form $(S,T]$ where $0\leq S<T$. We do
not consider this latter question in the article. We finally recall
that when $T<\mathcal{T}(u_d^{u_0})$ the LDP gives that
$$\lim_{\epsilon\rightarrow0}\epsilon\log\mathbb{P}\left(\mathcal{T}(u^{\epsilon,u_0})>
T\right)=0,$$ indeed this is not a large deviation event, we obtain
similarly that when $T>\mathcal{T}(u_d^{u_0})$ the LDP gives that
$$\lim_{\epsilon\rightarrow0}\epsilon\log\mathbb{P}\left(\mathcal{T}(u^{\epsilon,u_0})
\leq T\right)=0.$$

\subsection{Probability that blow-up occur before time T}\label{s52}
\begin{prpstn}\label{prpstn6}
If $T<\mathcal{T}_K^i=\inf_{u_0\in
K}\mathcal{T}\left(u_d^{u_0}\right)$, where $K\subset\subset\H1$,
then there exists $c$ positive such that
$$\overline{\lim}_{\epsilon\rightarrow0}\epsilon\log\sup_{u_0\in K}\mathbb{P}
\left(\mathcal{T}(u^{\epsilon,u_0})\leq T\right)\leq-c.$$
\end{prpstn}
\begin{proof}
Since $\mathcal{T}$ is lower semicontinuous,
$\mathcal{T}^{-1}([0,T])$ is a closed set. Suppose now that there
exists a sequence $(u_n,h_n)$ in $K\times\xLtwo(0,+\infty;\L2)$ such
that $$\mathcal{T}\left(S^c(u_n,h_n)\right)\leq T$$ and
$\lim_{n\rightarrow\infty} h_n=0$. Since $K$ is a compact set we may
extract a subsequence $u_{\varphi(n)}$ such that $u_{\varphi(n)}$
converges to some $\tilde{u}$. Also, if we denote by
$f_n(\cdot)=\int_0^{\cdot}\Phi h_n(s)ds$, $f_n$ converges to zero in
$\xCn\left([0,+\infty);\xHn^s(\R^d,\R)\right)$ and satisfies
$\mathcal{T}\left(S(u_n,f_n)\right)\leq T$. Also there exists $a$
positive such that for every $n$ in $\N$, $f_n$ in $C_a$. The
semicontinuity of $\mathcal{T}$ along with Proposition \ref{prpstn3}
give that
$$T\geq\underline{\lim}_{n\rightarrow\infty}\mathcal{T}
\left(S(u_{\varphi(n)},f_{\varphi(n)})\right)
\geq\mathcal{T}\left(S(\tilde{u},0)\right)\geq\mathcal{T}_K^i>T,$$
which is contradictory.
\end{proof}
When $K=\{u_0\}$ we obtain the same result as in Proposition 5.5 of
\cite{EG}.

\subsection{Probability that blow-up occur after time T}\label{s51}
In the following we consider the case $d=2$ or $d=3$ and a cubic
nonlinearity, i.e. $\sigma=1$. In that case blow-up may occur.
\begin{prpstn}\label{prpstn7} Let $U^{u_0}$ be the solution of the
free Schr\"odinger equation
with initial data $u_0$ in $\xHn^r(\R^d)$ where $r>\frac{d}{2}\vee
s$, assume that $T>\mathcal{T}\left(u_d^{u_0}\right)$ and that
$\spn\left\{|U^{u_0}(t)|^2,t\in[0,2T]\right\}$ belongs to the range
of $\Phi$. There exists $c$ positive such that
$$\underline{\lim}_{\epsilon\rightarrow0}\epsilon\log\mathbb{P}
\left(\mathcal{T}(u^{\epsilon,u_0})> T\right)\geq-c.$$
\end{prpstn}
\begin{rmrk}
It is known for example that for some Gaussian initial data $u_0$,
see for example \cite{BSDDM}, the solutions of NLS blows up in
finite time. Also, the solutions of the free equation are smooth and
strongly decreasing at infinity, thus it may easily be checked that
it is possible to define an Hilbert-Schmidt operator $\Phi$ such
that the last assumption holds.
\end{rmrk}
\begin{proof} Define $F^{u_0}$ by $F^{u_0}(t)=-\int_0^{t\wedge2T}
|U^{u_0}(s)|^{2}ds$.
The control is such that $S(u_0,F^{u_0})=U^{u_0}$ on $[0,2T]$, which
does not blow up, thus $\mathcal{T}\left(S(u_0,F^{u_0})\right)\geq
2T$. Also, $F^{u_0}$ belongs to
$\xCn\left([0,+\infty),\xHn^s(\R^d,\R)\right)$ since for
$r>\frac{d}{2}$ $\xHn^r(\R^d)$ is an algebra and $U^{u_0}$ belongs
to $\xCn([0,+\infty),\xHn^r(\R^d))$. Finally, from the assumption on
$\Phi$, there exists $h$ in $\xLtwo(0,+\infty;\xLtwo(\R^d))$ setting
$h=0$ after $2T$ such that $\Phi h(s)=|U^{u_0}(s)|^{2}\indic_{s\leq
2T}$ and $F^{u_0}$ belongs to $C_a$ for some $a$ positive. We thus
obtain that $F^{u_0}$ belongs to
$$\left\{f\in\xCn\left([0,+\infty),\xHn^s(\R^d,\R)\right):
\mathcal{T}\left(S(u_0,f)\right)>T\right\}$$
and that $I^W(F^{u_0})\leq a<+\infty$.
\end{proof}
\begin{rmrk}
We obtain a result on compact sets $K$ in $\xHn^r(\R^d)$ for
$T>\sup_{u_0\in K}\mathcal{T}\left(u_d^{u_0}\right)$ provided that
$\spn\left\{|U^{u_0}(t)|^2,t\in[0,2T],u_0\in K\right\}$ belongs to
the range of $\Phi$ restricted to a ball of
$\xLtwo\left(0,2T;\xLtwo(\R^d)\right)$.
\end{rmrk}

\end{document}